\documentclass{article}%
\usepackage{amsfonts}
\usepackage{amsmath}
\usepackage{amssymb}
\usepackage{amsxtra}
\usepackage{graphicx}
\usepackage{geometry}
\usepackage{color}
\usepackage{colortbl}
\usepackage{caption}
\usepackage[draft]{varioref}%
\providecommand{\U}[1]{\protect\rule{.1in}{.1in}}
\newtheorem{theorem}{Theorem}

\newtheorem{lemma}[theorem]{Lemma}

\newtheorem{proposition}[theorem]{Proposition}
\newtheorem{remark}[theorem]{Remark}

\numberwithin{equation}{section}
\begin{document}

\title{The localization effect for eigenfunctions of the mixed boundary value problem
in a thin cylinder with distorted ends}
\author{G.Cardone\\University of Sannio - Department of Engineering\\Corso Garibaldi, 107 - 82100 Benevento, Italy\\email: giuseppe.cardone@unisannio.it
\and T.Durante\\University of Salerno\\Department of Information
Engineering and Applied Mathematics\\Via Ponte don Melillo, 84084
Fisciano (SA), Italy
\\email: durante@diima.unisa.it
\and S.A.Nazarov\\Institute of Mechanical Engineering Problems\\V.O., Bolshoi pr., 61, 199178, St. Petersburg, Russia.\\email: srgnazarov@yahoo.co.uk}
\maketitle

\begin{abstract}
A simple sufficient condition on curved end of a straight cylinder
is found that provides a localization of the principal
eigenfunction of the mixed boundary value for the Laplace operator
with the Dirichlet conditions on the lateral side. Namely, the
eigenfunction concentrates in the vicinity of the ends and decays
exponentially in the interior. Similar effects are observed in the
Dirichlet and Neumann problems, too.

\medskip

\medskip

\end{abstract}

\section{Introduction}

\subsection{Formulation of the spectral problem}

\begin{figure}
[ptb]
\begin{center}
\includegraphics[
height=0.6002in,
width=3.5734in
]%
{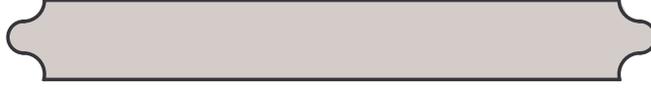}%
\caption{Thin cylinder with distorted ends}%
\label{f1}%
\end{center}
\end{figure}

Let $\omega\subset\mathbb{R}^{n-1}$ be a domain bounded by a simple closed
Lipschitz contour $\partial\omega$ and $n\geq2.$ Let also $H_{\pm}$ be
Lipschitz functions in $\overline{\omega}=\omega\cup\partial\omega.$ Given a
small parameter $h>0,$ we introduce the thin finite cylinder (Fig. \ref{f1})%
\begin{equation}
\Omega^{h}=\left\{  x=\left(  y,z\right)  \in\mathbb{R}^{n-1}\times
\mathbb{R}:\eta:=h^{-1}y\in\mathbb{\omega},\ \pm z<1\pm hH_{\pm}\left(
\eta\right)  \right\}  \label{1.1}%
\end{equation}
with the lateral side $\Sigma^{h}$ and the curved ends $\Gamma_{\pm}^{h},$%
\begin{equation}
\Sigma^{h}=\partial\Omega^{h}\diagdown\left(  \overline{\Gamma_{+}^{h}%
\cup\Gamma_{-}^{h}}\right)  ,\ \Gamma_{\pm}^{h}=\left\{  x:\eta\in\omega,\ \pm
z=1\pm hH_{\pm}\left(  \eta\right)  \right\}  . \label{1.2}%
\end{equation}
%

In the domain (\ref{1.1}) we consider the spectral mixed boundary value
problem%
\begin{align}
-\Delta_{x}u^{h}\left(  x\right)   &  =\lambda^{h}u^{h}\left(  x\right)
,\ \ x\in\Omega^{h},\label{1.3}\\
u^{h}\left(  x\right)   &  =0,\ \ x\in\Sigma^{h},\label{1.34}\\
\partial_{n}u^{h}\left(  x\right)   &  =0,\ \ x\in\Gamma_{\pm}^{h},
\label{1.4}%
\end{align}
where $\Delta_{x}=\Delta_{y}+\partial_{z}^{2}$ and $\Delta_{y}$ are the
Laplacians in $\mathbb{R}^{n}$ and $\mathbb{R}^{n-1},$ with $\partial
_{z}=\frac{\partial}{\partial z},$ and $\partial_{n}$ stands for
differentiation along the outward normal $n$ defined almost everywhere on the
Lipschitz surfaces $\Gamma_{\pm}^{h}.$ The variational formulation of the
problem (\ref{1.3})-(\ref{1.4}) reads: to find $\lambda^{h}\in\mathbb{R}$ and
$u^{h}\in\mathring{H}^{1}\left(  \Omega^{h};\Sigma^{h}\right)  ,$
$u^{h}\not \equiv 0,$ such that the integral identity \cite{Lad}%
\begin{equation}
\left(  \nabla_{x}u^{h},\nabla_{x}v\right)  _{\Omega^{h}}=\lambda\left(
u^{h},v\right)  _{\Omega^{h}},\ \ v\in\mathring{H}^{1}\left(  \Omega
^{h};\Sigma^{h}\right)  , \label{1.5}%
\end{equation}
is valid. Here $\left(  \ ,\ \right)  _{\Omega^{h}}$ is the natural inner
product in the Lebesgue space $L^{2}\left(  \Omega^{h}\right)  $ and
$\mathring{H}^{1}\left(  \Omega^{h};\Sigma^{h}\right)  $ the Sobolev space of
functions in $\Omega^{h}$ meeting the Dirichlet condition (\ref{1.34}).

The spectral boundary value problem (\ref{1.3})-(\ref{1.4}) is of
interest in many applicable disciplines related to acoustics and
electromagnetism. Since $h$ is a dimensionless parameter, the thin
domain (\ref{1.1}) of length $2+O\left(  h\right)  $ can also be
regarded as a long tubular domain with a
transverse cross-section of unit diameter. In this way problem (\ref{1.3}%
)-(\ref{1.4}) can be understood as the spectral problem for a cylindrical
waveguide with soft walls and hard ends (see \cite{W, LM} and others).

Since the left-hand side of (\ref{1.5}) serves as an inner product in
$\mathring{H}^{1}\left(  \Omega^{h};\Sigma^{h}\right)  $ and the embedding
$\mathring{H}^{1}\left(  \Omega^{h};\Sigma^{h}\right)  \subset L^{2}\left(
\Omega^{h}\right)  $ is compact, the spectral problem (\ref{1.5}) admits the
positive unbounded sequence of eigenvalues
\begin{equation}
0<\lambda_{1}^{h}<\lambda_{2}^{h}\leq\lambda_{3}^{h}\leq...\leq\lambda_{j}%
^{h}\leq...\rightarrow+\infty\label{1.6}%
\end{equation}
where the convention on repeated multiple eigenvalues is accepted. The
corresponding eigenfunctions $u_{1}^{h},u_{2}^{h},u_{3}^{h},...,u_{j}^{h},...$
in $\mathring{H}^{1}\left(  \Omega^{h};\Sigma^{h}\right)  $ can be subject to
the orthogonality and normalization conditions%
\begin{equation}
\left(  u_{j}^{h},u_{k}^{h}\right)  _{\Omega^{h}}=\delta_{j,k},\ \ j,k\in
\mathbb{N}=\left\{  1,2,3,...\right\}  , \label{1.7}%
\end{equation}
where $\delta_{j,k}$ is the Kronecker symbol. According to the strong maximum
principle, the first eigenvalue $\lambda_{1}^{h}$ is simple while the
corresponding eigenfunction $u_{1}^{h}$ may be chosen positive in $\Omega
^{h}.$

In the case $H_{\pm}=0$ the cylinder has straight ends and the explicit
dependence on $h$ can be clarified:%
\begin{equation}
\lambda_{p,q}^{h}=h^{-2}\mu_{p}+\frac{\pi^{2}q^{2}}{4},\ \ u_{p,q}^{h}\left(
x\right)  =\varphi_{p}\left(  y\right)  \cos\left(  \frac{\pi q}{2}\left(
z+1\right)  \right)  . \label{1.8}%
\end{equation}
The eigenpairs (\ref{1.8}) are renumerated with the two indices $p,q\in
\mathbb{N}.$ Furthermore, $\varphi_{1},\varphi_{2},\varphi_{3},...,\varphi
_{j},...$ are eigenfunctions corresponding to the eigenvalues%
\begin{equation}
0<\mu_{1}<\mu_{2}\leq\mu_{3}\leq...\leq\mu_{j}\leq...\rightarrow
+\infty\label{1.9}%
\end{equation}
of the spectral Dirichlet problem on the cross section%
\begin{equation}
-\Delta_{y}\varphi\left(  y\right)  =\mu\varphi\left(  y\right)  ,\ y\in
\omega,\ \ \ \varphi\left(  y\right)  =0,\ y\in\partial\omega, \label{1.10}%
\end{equation}
or, what is the same,
\begin{equation}
\left(  \nabla_{y}\varphi,\nabla_{y}\psi\right)  _{\omega}=\mu\left(
\varphi,\psi\right)  _{\omega},\ \ \psi\in\mathring{H}^{1}\left(
\omega\right)  :=\mathring{H}^{1}\left(  \omega;\partial\omega\right)  .
\label{1.11}%
\end{equation}
Similarly to (\ref{1.7}), the orthogonality and normalization conditions are
satisfied:%
\begin{equation}
\left(  \varphi_{j},\varphi_{k}\right)  _{\omega}=\delta_{j,k},\ j,k\in
\mathbb{N}. \label{1.12}%
\end{equation}

\subsection{The localization of eigenfunctions\label{sect1.2}}

All the eigenfunctions in (\ref{1.8}) are oscillating and do not become
infinitesimal as $h\rightarrow+\infty$ in any fragment $\omega_{h}%
\times\left(  -l_{-},l_{+}\right)  $ of the cylinder $\Omega_{h}=\omega
_{h}\times\left(  -1,1\right)  ;$ here $\omega_{h}=\left\{  y:\eta\in
\omega\right\}  $ and $l_{\pm}$ are arbitrarily fixed numbers, $-1\leq
-l_{-}<l_{+}\leq1.$ The main goal of the paper is to describe the
\textit{localization effect} for eigenfunction in the cylinder with distorted
ends. In other words, we reveal the profile functions $H_{\pm}$ in (\ref{1.1})
such that at least the first eigenfunction $u_{1}^{h}$ concentrates in the
vicinity of the ends and is of the exponential small order $\exp\left(
-h^{-1}\tau\right)  ,$ $\tau>0,$ outside a neighborhood of $\Gamma_{\pm}^{h}.$

\begin{figure}
[ptb]
\begin{center}
\includegraphics[
height=0.761in,
width=3.3996in
]%
{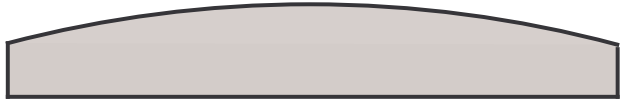}%
\caption{Thin curved trapezoid}%
\label{f2}%
\end{center}
\end{figure}

A similar, but of other kind, localization effect appears in the Dirichlet
problem for the Helmgoltz equation in a cylinder with the varying
cross-section (see, e.g., \cite{Sol12}). To outline this kind of the
localization, we briefly consider the Dirichlet problem in the thin curved
trapezoid (Fig. \ref{f2})%
\begin{equation}
\Omega^{h}=\left\{  x=\left(  y,z\right)  \in\mathbb{R}^{2}:\left\vert
z\right\vert <1,\ h^{-1}y\in\left(  0,H\left(  z\right)  \right)  \right\}
\label{1.20}%
\end{equation}
where $H\in C^{2}\left[  -1,1\right]  $ is a positive function with the strict
($b=-\partial_{z}^{2}H\left(  0\right)  >0$) global maximum at $z=0,$ monotone
for $z\in\left[  -1,0\right]  $ and $z\in\left[  0,1\right]  .$%

We present here only primary formal asymptotic analysis and refer to
\cite{na259, na310} for details and much more scrupulous explanation and to
\cite{Sol12} and \cite{BouPiat} for different approaches based on the spectral
theory of operators in Hilbert space and the $\Gamma$-convergence technique,
respectively. Accepting the asymptotic ans\"{a}tze for eigenpairs%
\begin{align}
\lambda^{h}  &  =h^{-2}\pi^{2}H\left(  0\right)  ^{-2}+h^{-1}\Lambda
+\widetilde{\lambda}^{h},\label{1.21}\\
u^{h}\left(  x\right)   &  =h^{-3/4}U\left(  h^{-1/2}z\right)  \sin\left(
h^{-1}\pi H\left(  0\right)  ^{-1}y\right)  +\widetilde{u}^{h}\left(
x\right)  ,\nonumber
\end{align}
we use the standard procedure for the dimension reduction in the rapid
variables $\eta=h^{-1}y,$ $\zeta=h^{-1/2}z.$ We insert the ans\"{a}tze
(\ref{1.21}) into the equation (\ref{1.3}) and the Neumann boundary conditions
(\ref{1.34}). Then we differentiate and apply the Taylor formula%
\[
H\left(  z\right)  ^{-2}=H\left(  0\right)  ^{-2}\left(  1+z^{2}H\left(
0\right)  ^{-1}b+O\left(  \left\vert z\right\vert ^{3}\right)  \right)
=H\left(  0\right)  ^{-2}\left(  1+h\zeta^{2}H\left(  0\right)  ^{-1}%
b+O\left(  h^{3}\left\vert \zeta\right\vert ^{3}\right)  \right)
\]
(see, e.g., \cite[\S 4]{na259} and \cite{na310} for details). Factors on
$h^{-2}$ in the differential equation cancell each other due to the first
formula (\ref{1.21}). Finally, factors on $h^{-1}$ mould the limit ordinary
differential equation%
\begin{equation}
-\partial_{\zeta}^{2}U\left(  \zeta\right)  +B\zeta^{2}U\left(  \zeta\right)
=\Lambda U\left(  \zeta\right)  ,\ \zeta\in\mathbb{R}, \label{1.22}%
\end{equation}
with the coefficient $B=H\left(  0\right)  ^{-3}b>0.$ Eigenpairs of the
spectral problem (\ref{1.22}), describing the harmonic oscillator, are known
(see, e.g., \cite{LL}):%
\begin{align}
\Lambda_{j}  &  =B^{1/2}\left(  2j+1\right)  ,\label{1.23}\\
U_{j}\left(  \zeta\right)   &  =\exp\left(  \frac{1}{2}B^{1/2}\zeta
^{2}\right)  \left(  \frac{d}{d\zeta}\right)  ^{j}\exp\left(  -B^{1/2}%
\zeta^{2}\right)  ,\ j\in\mathbb{R}.\nonumber
\end{align}
The formulas (\ref{1.23}) complete the asymptotic ans\"{a}tze which exhibit
the first asymptotic series of eigenvalues in the thin curved strip
(\ref{1.20}) while estimates for the remainders $\widetilde{\lambda}^{h}$ and
$\widetilde{u}^{h}$ in (\ref{1.21}) can be derived by different approaches
(see \cite{Sol12, na259, BouPiat} and others). Similar series of eigenvalues
with stable asymptotics may be derived, e.g., by using local maxima of the
function $H$ in (\ref{1.20}).

The distinguishing feature of the function $U_{j}$ in (\ref{1.23}) is the
exponential decay $o\left(  \exp\left(  -\sigma\zeta^{2}\right)  \right)  $
for $\zeta\rightarrow\pm\infty;$ here $\sigma\in\left(  0,\frac{1}{2}%
B^{1/2}\right)  .$ The corresponding eigenfunction $u_{j}^{h}$ concentrates in
the vicinity of the coordinate origin $\mathcal{O}$ and decays as $o\left(
\exp\left(  -h^{-1}\sigma z^{2}\right)  \right)  $ at a distance from
$\mathcal{O}.$ On the other hand, the first eigenfunction $u_{1}^{h}$ of the
problem (\ref{1.3})-(\ref{1.4}) in the rectangle $\left(  0,h\right)
\times\left(  -1,1\right)  $ oscillates and is spread over the whole domain.
In the trapezoid (\ref{1.20}) with the concave upper side (Fig. \ref{f2})
$u_{1}^{h}$ is localized in the $ch^{1/2}-$neighborhood of the highest point
in the upper side.

Localization effects of this kind have been described in the papers
\cite{na259, na310, Sol12, BouPiat, AlPi1, AlPi2} for eigenfunctions of
similar and other singular perturbed boundary value problems in domains in
$\mathbb{R}^{n}.$

In the present paper we study the localization effect of the other kind in the
$ch-$neighborhood of the ends $\Gamma_{\pm}^{h}$ of the cylinder $\Omega^{h}.$
The localization of eigenfunctions relies upon the boundary layer phenomenon
and the discrete spectrum of the limit boundary value problem (\ref{3.3}%
)-(\ref{3.5}) in the semi-infinite cylinders%
\begin{equation}
\Pi_{\pm}=\left\{  \xi^{\pm}=\left(  \eta^{\pm},\zeta^{\pm}\right)  :\eta
^{\pm}\in\omega,\ \zeta^{\pm}>-H_{\pm}\left(  \eta\right)  \right\}
\label{1.25}%
\end{equation}
obtained from the thin domain (\ref{1.1}) by the coordinate dilation
(\ref{3.1}) and the formal setting $h=0$. It is known (see, e.g., \cite{W})
can be readily verified that the continuous spectrum of the problem in
$\Pi_{\pm}$ coincides with the ray $\left[  \mu_{1},+\infty\right)  $ but the
discrete spectrum can appear below the cut-off $\mu_{1}$ which implies the
first eigenvalue of the problem (\ref{1.1}). Points $\Lambda^{\pm}$ of the
discrete spectrum in $\left(  0,\mu_{1}\right)  $ give rise to the so-called
\textit{trapped modes}, i.e. solutions of the homogeneous problem
(\ref{3.3})-(\ref{3.5}) with the exponential decay at infinity (see \cite{U}
and, e.g., reviews in \cite{LM, French, NaSLS}). These trapped modes become
the main asymptotic term in expansions of eigenfunctions of the problem
(\ref{1.5}). At the same time, the eigenvalues $\Lambda^{\pm}\in\left(
0,\mu_{1}\right)  $ of the problem (\ref{3.3})-(\ref{3.5}) in $\Pi_{\pm}$
after the multiplication by $h^{-2}$ approximate some eigenvalues in
(\ref{1.6}) with the precision $O\left(  \exp\left(  -h^{-1}\tau\right)
\right)  ,$ $\tau>0.$

We especially mention the paper \cite{na259} (see also \cite{na310}) where
both the approaches discussed above are combined to prove that the first
eigenfunction of the mixed boundary value problem for the equation (\ref{1.3})
in a thin plate $\Omega^{h}\subset\mathbb{R}^{3}$ with the distorted lateral
side may be localized near the point of the maximal curvature of the
longitudinal cross-section. Note that in this case the decay rate is of order
$\exp\left(  -h^{-1}\tau\right)  $ inside the domain, but $\exp\left(
-h^{-1/2}\sigma\right)  $ along the lateral side.

\subsection{Structure of the paper}

In Section \ref{sect2} after a brief comment on the continuous and discrete
spectra of the problem in the cylinder; we derive a simple sufficient
conditions for the existence of eigenvalues in the interval $\left(  0,\mu
_{1}\right)  $ below the continuous spectrum (see Theorems \ref{theor3.1} and
\ref{theor3.2}).

In Section \ref{sect3} we display asymptotic ans\"{a}tze for eigenvalues and
eigenfunctions of the problem (\ref{1.3})-(\ref{1.4}) and estimate the
asymptotic remainders (Theorem \ref{theor4.1}). Our proof is based on certain
estimates of Sobolev weighted norms of solutions to problems in $\Omega^{h}$
and $\Pi_{\pm}$ (Lemmas \ref{lem4.2} and \ref{lem5.2}) and differs from
approaches in \cite{na259, Sol12, BouPiat, AlPi1}. Also general results in the
theory of self-adjoint operators in Hilbert space are used throughout the paper.

\begin{figure}
[ptb]
\begin{center}
\includegraphics[
height=0.5068in,
width=2.7371in
]%
{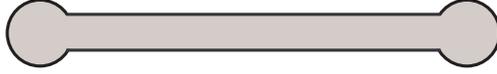}%
\caption{The dumbbell domain}%
\label{f3}%
\end{center}
\end{figure}

In Section \ref{sect4} we discuss similar spectral problems. First, we
consider the two-dimensional problem (\ref{1.3})-(\ref{1.4}), for which the
sufficient condition of nonempty discrete spectra only requires the negativity
of a coefficient in the Fourier series for the profile functions $\left(
0,1\right)  \ni\eta\mapsto H_{\pm}\left(  \eta\right)  $ (cf. formula
(\ref{7.2})). We also find out a domain where the first eigenfunction
concentrates at the both ends simultaneously. Then we consider the Dirichlet
and Neumann problems in $\Omega^{h}$ for the equation (\ref{1.3}). A
localization of eigenfunctions does not occur in the domain (\ref{1.1}) with
the Dirichlet conditions on the whole boundary (see Section \ref{sect4.2}).
However, we demonstrate the same localization effect in the dumbbell domain
(Fig. \ref{f3}). For the Neumann problem in the domain (\ref{1.1}), only
eigenfunction with large indices $j=O\left(  h^{-1}\right)  $ can admit the
localization in the vicinity of $\Gamma_{\pm}^{h}$ and we detect the
localization under a symmetry assumption on the domain $\Omega^{h}.$%

\section{The boundary layer phenomenon\label{sect2}}

\subsection{The mixed boundary value problem in a semi-cylinder.}

We introduce the stretched coordinates%
\begin{equation}
\xi^{\pm}=\left(  \eta,\zeta^{\pm}\right)  =\left(  h^{-1}y,h^{-1}\left(  1\mp
z\right)  \right)  . \label{3.1}%
\end{equation}
Since $\Delta_{x}=h^{-2}\Delta_{\xi^{\pm}},$ the change of the spectral
parameter%
\begin{equation}
\lambda^{h}\mapsto h^{-2}\Lambda\label{3.2}%
\end{equation}
and the formal setting $h=0$ in the original problem (\ref{1.3})-(\ref{1.4})
leads formally to the limit problem in a semi-infinite cylinder (\ref{1.25})
with a curvilinear end%
\begin{align}
-\Delta_{\xi^{\pm}}U  &  =\Lambda U\text{ \ in \ }\Pi_{\pm},\label{3.3}\\
U  &  =0\text{\ \ on \ }\Sigma_{\pm}=\partial\Pi_{\pm}\diagdown\overline
{\Gamma_{\pm}},\label{3.4}\\
\partial_{\nu}U  &  =0\text{ \ on \ }\Gamma_{\pm}=\left\{  \xi^{\pm}:\eta
\in\omega,\ \zeta^{\pm}=-H_{\pm}\left(  \eta\right)  \right\}  . \label{3.5}%
\end{align}
The variational formulation of the mixed boundary value problem reads: to find
$\Lambda\in\mathbb{R}_{+}$ and $U\in\mathring{H}^{1}\left(  \Pi;\Sigma\right)
,$ $U\not \equiv 0,$ such that%
\begin{equation}
\left(  \nabla_{\xi}U,\nabla_{\xi}V\right)  _{\Pi}=\Lambda\left(  U,V\right)
_{\Pi},\ \ V\in\mathring{H}^{1}\left(  \Pi;\Sigma\right)  . \label{3.99}%
\end{equation}
Here and in the sequel we omit the index $\pm$ in the notation.

We relate the spectral problem (\ref{3.3})-(\ref{3.5}) with the positive
self-adjoint unbounded operator $A$ in the Hilbert space $L^{2}\left(
\Pi\right)  $ with the domain $\mathcal{D}\left(  A\right)  $ generated,
according to \cite[\S 10.1, 10.2]{BiSo}, by the symmetric quadratic form%
\begin{equation}
a\left(  U,V\right)  =\left(  \nabla_{\xi}U,\nabla_{\xi}V\right)  _{\Pi
},\ \ U,V\in\mathring{H}^{1}\left(  \Pi;\Sigma\right)  . \label{3.6}%
\end{equation}

\begin{figure}
[ptb]
\begin{center}
\includegraphics[
height=0.7057in,
width=4.8066in
]%
{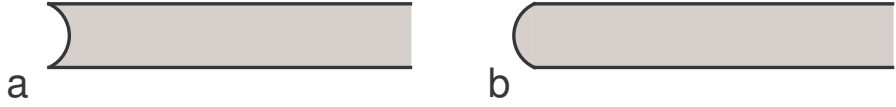}%
\caption{Acute and obtuse edges}%
\label{f4}%
\end{center}
\end{figure}

\begin{remark}
\label{rem3.11}If the contour $\partial\omega$ and the function $H_{\pm}$ are
smooth, the operator $A$ appears as the closure of the operator $A_{0}$ with
the differential expression $-\Delta_{\xi}$ and the domain $\mathcal{D}\left(
A_{0}\right)  =H^{2}\left(  \Pi\right)  \cap\mathring{H}^{1}\left(  \Pi
;\Sigma\right)  .$ From a result in \cite{na101} based on the theory of
elliptic problems in domains with piecewise smooth boundaries (see, e.g.,
\cite{NaPl, KoMaRo} and, particularly, \cite{Ko, Ko1}), shows that the domain
$\mathcal{D}\left(  A\right)  $ coincides with $\mathcal{D}\left(
A_{0}\right)  $ in the case $\partial_{n}H_{\pm}>0$ on $\partial\omega,$ i.e.
for the acute edge $\gamma_{\pm}=\left\{  \xi:\eta\in\partial\omega
,\ \zeta^{\pm}=-H_{\pm}\left(  \eta\right)  \right\}  $ on the surface
$\partial\Pi$ (cf. Fig. \ref{f4},a). If $\partial_{n}H_{\pm}<0$ and,
therefore, the edge $\gamma_{\pm}$ is obtuse (cf. Fig. \ref{f4},b), then
$\mathcal{D}\left(  A_{0}\right)  \subsetneqq\mathcal{D}\left(  A\right)  .$
For any edge, $\mathcal{D}\left(  A\right)  $ falls into the Kondratiev space
\cite{Ko}
\[
\left\{  U\in\mathring{H}^{1}\left(  \Pi_{\pm};\Gamma_{\pm}\right)  :\left(
1+\text{dist}\left(  \xi,\gamma_{\pm}\right)  \right)  ^{k-1}\nabla_{\xi}%
^{k}U\in L^{2}\left(  \Pi\right)  ,\ k=0,1,2\right\}  .\ \blacksquare
\]

\end{remark}%

\subsection{The continuous spectrum}

\begin{figure}
[ptb]
\begin{center}
\includegraphics[
height=1.5065in,
width=4.6406in
]%
{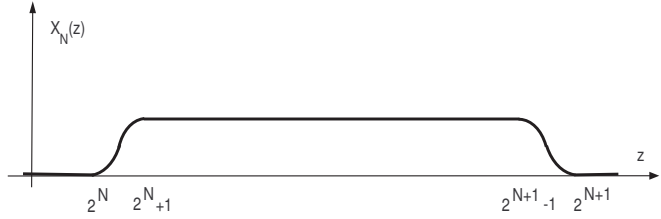}%
\caption{The plateau function}%
\label{f5}%
\end{center}
\end{figure}

Let us list properties of the spectrum $\sigma\left(  A\right)  $ of the
operator $A,$ which are well known (see, e.g., the review \cite{LM}) and
readily follow from general results on semi-bounded self-adjoint operators in
Hilbert space (see, e.g., \cite{BiSo}). The set $\sigma\left(  A\right)  $
belongs to the positive real semi-axis $\mathbb{R}_{+}$ in the complex plane
$\mathbb{C}$. The operator $A$ has the continuous spectrum%
\begin{equation}
\sigma_{c}\left(  A\right)  =\left[  \mu_{1},+\infty\right)  \label{3.7}%
\end{equation}
where $\mu_{1}>0$ is the first eigenvalue of the Dirichlet problem
(\ref{1.11}) on the cross-section $\omega$ (see also formulas (\ref{1.10}) and
(\ref{1.9})). To verify this fact, one may observe that, for $\Lambda\geq
\mu_{1},$ the function%
\begin{equation}
\left(  \eta,\zeta\right)  \mapsto\exp\left(  \pm i\left(  \Lambda-\mu
_{1}\right)  ^{1/2}\zeta\right)  \varphi_{1}\left(  \eta\right)  , \label{3.8}%
\end{equation}
where $\varphi_{1}$ is the eigenfunction corresponding to $\mu_{1}$ and $i$
the imaginary unit, satisfies the equation (\ref{3.3}) in the entire cylinder
$\omega\times\mathbb{R}$ and the Dirichlet conditions on the cylindrical
surface $\partial\omega\times\mathbb{R}$ but has no decay for both
$\zeta\rightarrow\pm\infty.$ Now multiplying (\ref{3.8}) by the plateau
functions $X_{N}$ with the graph in Fig. \ref{f5} provides the singular Weyl
sequence for the operator $A$ at the point $\Lambda$ (see, e.g.,
\cite[\S 9.1]{BiSo}, \cite[Thm. 3.1.1]{NaPl}, \cite{LM, NaSLS} and others).
This ensures that $\Lambda$ belongs to the essential spectrum $\sigma
_{e}\left(  A\right)  .$ The kernel of the operator $A-\Lambda$ is always
finite dimensional (see, e.g., \cite[Remark 3.1.5]{NaPl}) and, therefore,
$\Lambda\in\sigma_{c}\left(  A\right)  .$%

\subsection{The discrete spectrum\label{sect2.3}}

The interval $\left(  0,\mu_{1}\right)  \subset\mathbb{R}_{+}$ contains the
discrete spectrum $\sigma_{d}\left(  A\right)  $ only. According to \cite[Thm.
10.2.1]{BiSo}, the lower bound of the spectrum $\sigma\left(  A\right)  $ can
be computed as follows:%
\begin{equation}
\min\left\{  \Lambda:\Lambda\in\sigma\left(  A\right)  \right\}  =\inf
_{u\in\mathring{H}^{1}\left(  \Pi;\Sigma\right)  \diagdown\left\{  0\right\}
}\frac{a\left(  U,U\right)  }{\left\Vert U;L^{2}\left(  \Pi\right)
\right\Vert ^{2}}. \label{3.9}%
\end{equation}
We emphasize that the infimum on the left is taken over all nontrivial
functions in the domain $\mathring{H}^{1}\left(  \Pi;\Sigma\right)  $ of the
form (\ref{3.6}) which is bigger than the domain $\mathcal{D}\left(  A\right)
$ of the operator $A$ (cf. Remark \ref{rem3.11}). It clearly is one advantage
of the theory \cite[Ch. 10]{BiSo} applied here.

Our immediate objective becomes to find out a condition on the shape of the
end $\Gamma$ which provides a trial function $W\in\mathring{H}^{1}\left(
\Pi;\Sigma\right)  $ such that
\begin{equation}
a\left(  W,W\right)  <\mu_{1}\left\Vert W;L^{2}\left(  \Pi\right)  \right\Vert
^{2}. \label{3.10}%
\end{equation}
Indeed, the inequality (\ref{3.10}) ensures that the lower bound (\ref{3.9})
is less than the cut-off $\mu_{1}$ for the continuous spectrum $\sigma
_{c}\left(  A\right)  $ and, hence, the discrete spectrum $\sigma_{d}\left(
A\right)  $ cannot be empty. We employ an approach in \cite{na259, na322} (see
also \cite[\S 3.6]{NaSLS}). For any $\varepsilon>0,$ we set%
\begin{equation}
W\left(  \xi\right)  =\exp\left(  -\varepsilon\zeta\right)  \varphi_{1}\left(
\eta\right)  . \label{3.11}%
\end{equation}
Since $\overline{\Sigma}=\partial\Pi\cap\left(  \partial\omega\times
\mathbb{R}\right)  $ and $\varphi_{1}\in\mathring{H}^{1}\left(  \omega
;\partial\omega\right)  ,$ the function $W$ belongs to $\mathring{H}%
^{1}\left(  \Pi;\Sigma\right)  .$ We obtain%
\begin{align}
\left\Vert W;L^{2}\left(  \Pi\right)  \right\Vert ^{2}  &  =\int_{\omega}%
\int_{-H\left(  \eta\right)  }^{+\infty}\exp\left(  -2\varepsilon\zeta\right)
\varphi_{1}\left(  \eta\right)  ^{2}d\zeta d\eta=\nonumber\\
&  =\frac{1}{2\varepsilon}\int_{\omega}\exp\left(  2\varepsilon H\left(
\eta\right)  \right)  \varphi_{1}\left(  \eta\right)  ^{2}d\eta=\label{3.12}\\
&  =\frac{1}{2\varepsilon}+\int_{\omega}H\left(  \eta\right)  \varphi
_{1}\left(  \eta\right)  ^{2}d\eta+\varepsilon\int_{\omega}H\left(
\eta\right)  ^{2}\varphi_{1}\left(  \eta\right)  ^{2}d\eta+O\left(
\varepsilon^{2}\right)  ,\nonumber
\end{align}%
\begin{align*}
\left\Vert \nabla_{\xi}W;L^{2}\left(  \Pi\right)  \right\Vert ^{2}  &
=\int_{\omega}\int_{-H\left(  \eta\right)  }^{+\infty}\exp\left(
-2\varepsilon\zeta\right)  \left(  \left\vert \nabla_{\eta}\varphi_{1}\left(
\eta\right)  \right\vert ^{2}+\varepsilon^{2}\left\vert \varphi_{1}\left(
\eta\right)  \right\vert ^{2}\right)  d\zeta d\eta=\\
&  =\int_{\omega}\exp\left(  2\varepsilon H\left(  \eta\right)  \right)
\varphi_{1}\left(  \eta\right)  ^{2}d\eta=\\
&  =\frac{\mu_{1}}{2\varepsilon}+\int_{\omega}H\left(  \eta\right)  \left\vert
\nabla_{\eta}\varphi_{1}\left(  \eta\right)  \right\vert ^{2}d\eta
+\varepsilon\left(  \frac{1}{2}+\int_{\omega}H\left(  \eta\right)
^{2}\left\vert \nabla_{\eta}\varphi_{1}\left(  \eta\right)  \right\vert
^{2}d\eta\right)  +O\left(  \varepsilon^{2}\right)  .
\end{align*}
Here we have used the formulas%
\begin{equation}
\left\Vert \varphi_{1};L^{2}\left(  \omega\right)  \right\Vert
=1,\ \ \left\Vert \nabla_{\eta}\varphi_{1};L^{2}\left(  \omega\right)
\right\Vert =\sqrt{\mu_{1}}, \label{3.00}%
\end{equation}
which directly follow from (\ref{1.12}) and (\ref{1.11}).

Inserting (\ref{3.12}) into (\ref{3.10}), we see that the terms of order
$\varepsilon^{-1}$ cancel each other. Collecting terms of order $\varepsilon
^{0},$ we observe that the inequality (\ref{3.10}) is valid with the function
(\ref{3.11}) and a small $\varepsilon>0$ in the case%
\begin{equation}
\int_{\omega}H\left(  \eta\right)  \left(  \left\vert \nabla_{\eta}\varphi
_{1}\left(  \eta\right)  \right\vert ^{2}-\mu_{1}\varphi_{1}\left(
\eta\right)  ^{2}\right)  d\eta<0. \label{3.13}%
\end{equation}
In other words, if one succeeds to find out a function $H$ such that the
relation (\ref{3.13}) is met, one readily detects an eigenvalue $\Lambda
_{1}\in\left(  0,\mu_{1}\right)  $ of the operator $A$ and simultaneously of
the problem (\ref{3.3})-(\ref{3.5}).

\begin{figure}
[ptb]
\begin{center}
\includegraphics[
height=0.8449in,
width=5.8721in
]%
{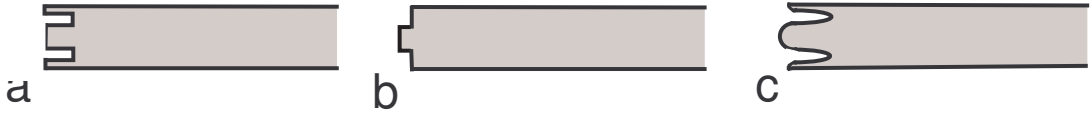}%
\caption{Shapes which provide trapped modes}%
\label{f6}%
\end{center}
\end{figure}

The function $\Phi=\left\vert \nabla_{\eta}\varphi_{1}\right\vert ^{2}-\mu
_{1}\left\vert \varphi_{1}\right\vert ^{2}$ is of mean zero (see (\ref{3.00}))
and, therefore, changes sign in $\omega.$ Hence, the inequality (\ref{3.13})
surely can be satisfied by an appropriate choice of the profile function $H.$
In this way $H$ can be of rather arbitrary behavior (see Fig. \ref{f6}).%

Note that $\Phi$ is positive in a neighborhood of the smooth contour
$\partial\omega$ because on one hand $\varphi_{1}=0$ on $\partial\omega$ and
on the other hand $\left\vert \nabla_{\eta}\varphi_{1}\right\vert >0$ by the
strong maximum principle. In Fig. \ref{f6},a we locate a \textit{negative} $H$
in this neighborhood. The set $\left\{  y\in\omega:\Phi\left(  y\right)
<0\right\}  $ is not empty and includes the support of a positive $H$ in Fig.
\ref{f6},b. Finally, $H$ is smooth and changes sign in Fig. \ref{f6},c.

\begin{theorem}
\label{theor3.1} If the function $H\in C\left(  \overline{\omega}\right)  $
meets the condition (\ref{3.13}) where $\left\{  \mu_{1},\varphi_{1}\right\}
$ is the first eigenpair of the Dirichlet problem (\ref{1.10}), the mixed
boundary value problem (\ref{3.99}) in the semi-infinite cylinder%
\begin{equation}
\Pi=\left\{  \xi=\left(  \eta,\zeta\right)  :\eta\in\omega,\ \zeta>-H\left(
\eta\right)  \right\}  \label{3.98}%
\end{equation}
has an eigenvalue below the cut-off $\mu_{1}$ for the continuous spectrum
(\ref{3.7}).
\end{theorem}

\subsection{A simplified sufficient condition for the existence of a trapped
mode}

Let us simplify the integral in (\ref{3.13}) under the smoothness assumption
$H\in C^{2}\left(  \overline{\omega}\right)  .$ Setting $\psi=H\varphi_{1}$ in
the integral identity (\ref{1.11}) for the spectral pair $\left\{  \mu
_{1},\varphi_{1}\right\}  $ yields%
\begin{align}
\mu_{1}\left(  \varphi_{1},H\varphi_{1}\right)  _{\omega}  &  =\left(
\nabla_{\eta}\varphi_{1},\nabla_{\eta}\left(  H\varphi_{1}\right)  \right)
_{\omega}=\left(  \nabla_{\eta}\varphi_{1},H\nabla_{\eta}\varphi_{1}\right)
_{\omega}+\left(  \nabla_{\eta}\varphi_{1},\varphi_{1}\nabla_{\eta}H\right)
_{\omega}=\label{3.97}\\
&  =\left(  \nabla_{\eta}\varphi_{1},H\nabla_{\eta}\varphi_{1}\right)
_{\omega}-\frac{1}{2}\left(  \varphi_{1}^{2},\Delta_{\eta}H\right)  _{\omega
}.\nonumber
\end{align}
Thus, the condition (\ref{3.13}) is equivalent to%
\begin{equation}
\int_{\omega}\varphi_{1}\left(  \eta\right)  ^{2}\Delta_{\eta}H\left(
\eta\right)  d\eta<0. \label{3.15}%
\end{equation}

\begin{theorem}
\label{theor3.2}If the function $H\in C^{2}\left(  \overline{\omega}\right)  $
meets the condition (\ref{3.15}) where $\varphi_{1}$ is the first
eigenfunction of the Dirichlet problem (\ref{1.10}), the mixed boundary value
problem (\ref{3.3})-(\ref{3.5}) in the cylinder (\ref{3.98}) has an eigenvalue
below the cut-off $\mu_{1}$ for the continuous spectrum (\ref{3.7}).
\end{theorem}

\begin{figure}
[ptb]
\begin{center}
\includegraphics[
height=0.8769in,
width=5.6507in
]%
{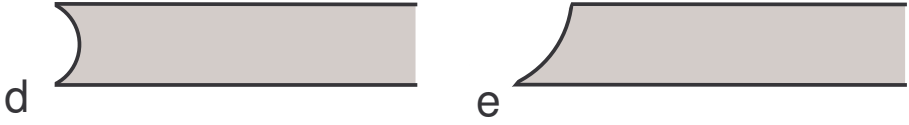}%
\caption{The superharmonic and subharmonic ends of the semi-infinite cylinder}%
\label{f7}%
\end{center}
\end{figure}

The inequality (\ref{3.15}) becomes true for a subharmonic function $H$ but
false for superharmonic. Fig. \ref{f7}, d and e, present two semi-cylinders
(\ref{1.25}) determined with a subharmonic function $H,$ while the condition
(\ref{3.13}) is not so evident for application in these cases. At the same
time, the semi-cylinder in Fig. \ref{f6}, c, is given by a function $H\in
C^{2}\left(  \overline{\omega}\right)  $ which is not subharmonic but an
eigenvalue $\Lambda\in\left(  0,\mu_{1}\right)  $ exists.

If $H$ is a harmonics, e.g., a linear function, the sufficient condition
(\ref{3.15}) is helpless. We attempt to make further use of the calculation
(\ref{3.12}) where we collect terms of order $\varepsilon.$ As a result, we
see that (\ref{3.10}) holds true with the function (\ref{3.11}) and a small
$\varepsilon>0$ provided%
\begin{equation}
\frac{1}{2}+\int_{\omega}H\left(  \eta\right)  ^{2}\left\vert \nabla_{\eta
}\varphi_{1}\left(  \eta\right)  \right\vert ^{2}d\eta-\mu_{1}\int_{\omega
}H\left(  \eta\right)  ^{2}\varphi_{1}\left(  \eta\right)  ^{2}d\eta<0.
\label{3.16}%
\end{equation}
We now repeat the calculation (\ref{3.97}) while changing $H$ for $H^{2}$ and
observe that the square of a harmonic function is superharmonic:%
\[
\frac{1}{2}\Delta_{\eta}H\left(  \eta\right)  ^{2}=H\left(  \eta\right)
\Delta_{\eta}H\left(  \eta\right)  +\left\vert \nabla_{\eta}H\left(
\eta\right)  \right\vert ^{2}=\left\vert \nabla_{\eta}H\left(  \eta\right)
\right\vert ^{2}\geq0.
\]
Hence, (\ref{3.16}) converts into the false inequality%
\[
\frac{1}{2}+\int_{\omega}\left\vert \nabla_{\eta}H\left(  \eta\right)
\right\vert ^{2}\varphi_{1}\left(  \eta\right)  ^{2}d\eta<0
\]
and, therefore, a harmonic function $H$ cannot assure the inequality
(\ref{3.10}) by the trial function (\ref{3.11}).%

\section{Justification of the asymptotics\label{sect3}}

\subsection{The hypothesis and the theorem on asymptotics of eigenvalues}

We assume that the spectral problems (\ref{3.3})-(\ref{3.5}) in the
semi-infinite cylinders $\Pi^{\pm}$ admit the eigenvalues%
\begin{equation}
\Lambda_{1}^{\pm},...,\Lambda_{N^{\pm}}^{\pm}\in\left(  0,\mu_{1}-\delta
_{1}\right]  . \label{4.1}%
\end{equation}
Theorems \ref{theor3.1} and \ref{theor3.2} provide sufficient conditions for
the existence of the discrete spectrum below the cut-off $\mu_{1}$ and,
accepting for $H_{\pm}$ either (\ref{3.15}), or (\ref{3.13}), we fix
$\delta_{1}>0$ such that the inclusion (\ref{4.1}) is valid and $N=N^{+}\cup
N^{-}\geq1.$ In this section we prove the following theorem.

\begin{theorem}
\label{theor4.1}If the hypothesis (\ref{4.1}) is true, the entries of the
eigenvalue sequence (\ref{1.6}) of the problem (\ref{1.5}) in the thin
cylinder (\ref{1.1}) satisfy the relations%
\begin{equation}
\left\vert \lambda_{p}^{h}-h^{-2}\Lambda^{\left(  p\right)  }\right\vert \leq
c_{p}\exp\left(  -\tau_{p}h^{-1}\right)  ,\text{ \ }h\in\left(  0,h_{p}%
\right]  , \label{4.2}%
\end{equation}
where $c_{p},$ $\tau_{p}$ and $h_{p}$ are certain positive numbers,
$p=1,...,N,$ and the eigenvalues (\ref{4.1}) are renumerated as follows:%
\begin{equation}
0<\Lambda^{\left(  1\right)  }\leq\Lambda^{\left(  2\right)  }\leq
...\leq\Lambda^{\left(  N\right)  }\leq\mu_{1}-\delta_{1}. \label{4.3}%
\end{equation}

\end{theorem}

An information on the corresponding eigenfunction $u_{1}^{h},...,u_{N}^{h}$ is
obtained in Section \ref{sect3.5} as well. The eigenfunctions $U_{j}^{\pm}$ of
the problems (\ref{3.99}) in $\Pi_{\pm}$ can be subject to the normalization
and orthogonality conditions%
\begin{equation}
\left(  U_{j}^{\pm},U_{k}^{\pm}\right)  _{\Pi_{\pm}}=\delta_{j,k}.
\label{4.312}%
\end{equation}
After the renumeration of the eigenvalues $\Lambda_{j}^{\pm}$ in (\ref{4.3}),
the corresponding eigenfunctions $U^{\left(  p\right)  },$ of course, keep
certain normalization and orthogonality conditions in $\Pi_{+}\times\Pi_{-}$
and we further still refer to (\ref{4.312}), although a new way of writing the
selfsame condition is accepted.

\subsection{A result on convergence}

We proceed with a primitive result on the convergence of the normalized
eigenvalues in the lower range of the spectrum (\ref{1.6}) and the
corresponding eigenfunctions.

Let $u^{h}\in\mathring{H}^{1}\left(  \Omega^{h};\Sigma^{h}\right)  $ be a
solution of the problem (\ref{1.5}) with
\begin{equation}
\lambda^{h}\in\left(  0,h^{-2}\left(  \mu_{1}-\tau_{0}^{2}\right)  \right]
,\ \tau_{0}>0. \label{4.4}%
\end{equation}
Aiming to estimate Sobolev weighted norm of $u^{h},$ we introduce the positive
weight function
\begin{equation}
\mathcal{R}_{\tau}\left(  x\right)  =\left\{
\begin{array}
[c]{c}%
\exp\left(  \tau h^{-1}\left(  1-\left\vert z\right\vert -c_{H}h\right)
\right)  ,\ x\in\Omega_{H}^{h},\\
1,\ \ \ \ \ \ \ \ \ \ \ \ \ \ \ \ \ \ \ \ \ \ x\in\Omega^{h}\diagdown
\Omega_{H}^{h},
\end{array}
\right.  \label{4.5}%
\end{equation}
where $\tau>0,$ $\Omega_{H}^{h}=\left\{  x=\left(  y,z\right)  :\eta\in
\omega,\ \left\vert z\right\vert <1-c_{H}h\right\}  $ and $c_{H}\geq
\max\left\vert H_{\pm}\left(  \eta\right)  \right\vert $ is fixed such that
$\Omega_{H}^{h}\subset\Omega^{h}.$ The function (\ref{4.5}) is continuous and
piecewise smooth while%
\begin{equation}
\left\vert \nabla_{x}\mathcal{R}_{\tau}\left(  x\right)  \right\vert \leq\tau
h^{-1}\mathcal{R}_{\tau}\left(  x\right)  ,\ \ x\in\Omega^{h},\ \ \nabla
_{x}\mathcal{R}_{\tau}\left(  x\right)  =0,\ \ x\in\Omega^{h}\diagdown
\Omega_{H}^{h}. \label{4.6}%
\end{equation}
We insert the test function $v=\mathcal{R}_{\tau}^{2}u\in\mathring{H}%
^{1}\left(  \Omega^{h};\Sigma^{h}\right)  $ into the integral identity
(\ref{1.5}) and obtain%
\begin{align}
\lambda^{h}\left\Vert \mathcal{R}_{\tau}u^{h};L^{2}\left(  \Omega^{h}\right)
\right\Vert ^{2}  &  =\lambda^{h}\left(  u^{h},\mathcal{R}_{\tau}^{2}%
u^{h}\right)  _{\Omega^{h}}=\left(  \nabla_{x}u^{h},\nabla_{x}\left(
\mathcal{R}_{\tau}^{2}u^{h}\right)  \right)  _{\Omega^{h}}=\label{4.7}\\
&  =\left(  \mathcal{R}_{\tau}\nabla_{x}u^{h},\nabla_{x}\left(  \mathcal{R}%
_{\tau}u^{h}\right)  \right)  _{\Omega^{h}}+\left(  \mathcal{R}_{\tau}%
\nabla_{x}u^{h},u^{h}\nabla_{x}\mathcal{R}_{\tau}\right)  _{\Omega_{H}^{h}%
}=\nonumber\\
&  =\left\Vert \nabla_{x}\left(  \mathcal{R}_{\tau}u^{h}\right)  ;L^{2}\left(
\Omega^{h}\right)  \right\Vert ^{2}-\left(  u^{h}\nabla_{x}\mathcal{R}_{\tau
},\nabla_{x}\left(  \mathcal{R}_{\tau}u^{h}\right)  \right)  _{\Omega_{H}^{h}%
}+\nonumber\\
&  +\left(  \nabla_{x}\left(  \mathcal{R}_{\tau}u^{h}\right)  ,u^{h}\nabla
_{x}\mathcal{R}_{\tau}\right)  _{\Omega_{H}^{h}}-\left\Vert u^{h}\nabla
_{x}\mathcal{R}_{\tau};L^{2}\left(  \Omega_{H}^{h}\right)  \right\Vert
^{2}.\nonumber
\end{align}
Two terms on the right-hand side cancel each other. We integrate over
$z\in\left(  -1+c_{H}h,1-c_{H}h\right)  $ the Friedrichs inequality in the
small domain $\omega^{h}=\left\{  y:h^{-1}y\in\omega\right\}  $ to observe
that%
\begin{equation}
\left\Vert \nabla_{y}w;L^{2}\left(  \Omega_{H}^{h}\right)  \right\Vert
^{2}\geq h^{-2}\mu_{1}\left\Vert w;L^{2}\left(  \Omega_{H}^{h}\right)
\right\Vert ^{2},\ \ w\in\mathring{H}^{1}\left(  \Omega^{h};\Sigma^{h}\right)
. \label{4.8}%
\end{equation}
Since $\mathcal{R}=1$ on $x\in\Omega^{h}\diagdown\Omega_{H}^{h}$ (see
(\ref{4.5})), from (\ref{4.6})-(\ref{4.8}) it follows that%
\begin{align}
\lambda^{h}\left\Vert u^{h};L^{2}\left(  \Omega^{h}\right)  \right\Vert ^{2}
&  \geq\lambda^{h}\left\Vert u^{h};L^{2}\left(  \Omega^{h}\diagdown\Omega
_{H}^{h}\right)  \right\Vert ^{2}=\lambda^{h}\left\Vert \mathcal{R}_{\tau
}u^{h};L^{2}\left(  \Omega^{h}\diagdown\Omega_{H}^{h}\right)  \right\Vert
^{2}=\label{4.888}\\
&  =\left\Vert \nabla_{x}\left(  \mathcal{R}_{\tau}u^{h}\right)  ;L^{2}\left(
\Omega^{h}\right)  \right\Vert ^{2}-\left\Vert u^{h}\nabla_{x}\mathcal{R}%
_{\tau};L^{2}\left(  \Omega_{H}^{h}\right)  \right\Vert ^{2}-\lambda
^{h}\left\Vert \mathcal{R}_{\tau}u^{h};L^{2}\left(  \Omega_{H}^{h}\right)
\right\Vert ^{2}\geq\nonumber\\
&  \geq h^{-2}\left(  \mu_{1}-\tau^{2}-h^{2}\lambda^{h}\right)  \left\Vert
\mathcal{R}_{\tau}u^{h};L^{2}\left(  \Omega_{H}^{h}\right)  \right\Vert
^{2}.\nonumber
\end{align}
By virtue of the assumption (\ref{4.4}), the factor on the last norm is bigger
than $h^{-2}\left(  \tau_{0}^{2}-\tau^{2}\right)  $ and stays positive in the
case $\tau\in\left(  0,\tau_{0}\right)  .$ Thus, reading the relation
(\ref{4.888}) without its middle part, we see that%
\begin{align}
\left\Vert \mathcal{R}_{\tau}u^{h};L^{2}\left(  \Omega_{H}^{h}\right)
\right\Vert ^{2}  &  \leq h^{2}\lambda^{h}\left(  \tau_{0}^{2}-\tau
^{2}\right)  ^{-1}\left\Vert u^{h};L^{2}\left(  \Omega^{h}\right)  \right\Vert
^{2}\leq\label{4.9}\\
&  \leq\mu_{1}\left(  \tau_{0}^{2}-\tau^{2}\right)  ^{-1}\left\Vert
u^{h};L^{2}\left(  \Omega^{h}\right)  \right\Vert ^{2}.\nonumber
\end{align}
Returning back to (\ref{4.7}), we use (\ref{4.9}) to conclude that%
\begin{align}
\left\Vert \nabla_{x}\left(  \mathcal{R}_{\tau}u^{h}\right)  ;L^{2}\left(
\Omega^{h}\right)  \right\Vert ^{2}  &  \leq\left(  \lambda^{h}+\tau^{2}%
h^{-2}\right)  \left\Vert \mathcal{R}_{\tau}u^{h};L^{2}\left(  \Omega
^{h}\right)  \right\Vert ^{2}\leq\label{4.999}\\
&  \leq h^{-2}\left(  \mu_{1}+\tau_{0}^{2}\right)  \left(  1+\mu_{1}\left(
\tau_{0}^{2}-\tau^{2}\right)  ^{-1}\right)  \left\Vert u^{h};L^{2}\left(
\Omega^{h}\right)  \right\Vert ^{2}.\nonumber
\end{align}
Taking into account that, by (\ref{4.6}) and (\ref{4.999}),
\[
\left\Vert \mathcal{R}_{\tau}\nabla_{x}u^{h};L^{2}\left(  \Omega^{h}\right)
\right\Vert ^{2}\leq c\left(  \left\Vert \nabla_{x}\left(  \mathcal{R}_{\tau
}u^{h}\right)  ;L^{2}\left(  \Omega^{h}\right)  \right\Vert ^{2}%
+h^{-2}\left\Vert \mathcal{R}_{\tau}u^{h};L^{2}\left(  \Omega^{h}%
\diagdown\Omega_{H}^{h}\right)  \right\Vert ^{2}\right)  \leq ch^{-2}%
\left\Vert \mathcal{R}_{\tau}u^{h};L^{2}\left(  \Omega^{h}\right)  \right\Vert
^{2},
\]
we formulate the obtained result.

\begin{lemma}
\label{lem4.2}Let $\left\{  \lambda^{h},u^{h}\right\}  $ be an eigenpair of
the problem (\ref{1.5}) such that $\left\Vert u^{h};L^{2}\left(  \Omega
^{h}\right)  \right\Vert =1$ and the inclusion (\ref{4.4}) holds true. Then
the following estimate is valid:%
\begin{equation}
h^{2}\left\Vert \mathcal{R}_{\tau}\nabla_{x}u^{h};L^{2}\left(  \Omega
^{h}\right)  \right\Vert ^{2}+\left\Vert \mathcal{R}_{\tau}u^{h};L^{2}\left(
\Omega^{h}\right)  \right\Vert ^{2}\leq c\left(  \tau\right)  , \label{4.10}%
\end{equation}
where $\mathcal{R}_{\tau}$ is the weight function (\ref{4.5}) with the
parameter $\tau\in\left(  0,\tau_{0}\right)  $ and the bound $c\left(
\tau\right)  $ is independent of $h\in\left(  0,1\right]  $ and $\left\{
\lambda^{h},u^{h}\right\}  $ but $c\left(  \tau\right)  \rightarrow+\infty$ as
$\tau\rightarrow\tau_{0}-0.$
\end{lemma}

Let the entry $\lambda_{j}^{h}$ of the eigenvalue sequence (\ref{1.6}) satisfy
the condition (\ref{4.4}). From the corresponding eigenfunction $u_{j}^{h}%
\in\mathring{H}^{1}\left(  \Omega^{h};\Sigma^{h}\right)  $ we construct the
two following functions in the semi-cylinders $\Pi_{\pm}\ni\xi^{\pm}:$%
\begin{equation}
U_{j\pm}^{h}\left(  \xi^{\pm}\right)  =h^{3/2}\chi\left(  h\zeta^{\pm}\right)
u_{j}^{h}\left(  h\eta,\mp h\zeta^{\pm}\pm1\right)  . \label{4.11}%
\end{equation}
Here $\xi^{\pm}$ are the stretched coordinates (\ref{3.1}) and $\chi\in
C^{\infty}\left(  \mathbb{R}\right)  $ is a cut-off function such that
$\chi\left(  t\right)  =1$ for $t\leq\frac{1}{3}$ and $\chi\left(  t\right)
=0$ for $t\geq\frac{2}{3},$ $0\leq\chi\leq1.$ Note that $h^{3/2}$ is a
normalization factor caused by the coordinate dilation $x\mapsto\xi$ that
together with Lemma \ref{lem4.2} furnish the relation%
\begin{align*}
\left\Vert U_{j+}^{h};L^{2}\left(  \Pi_{+}\right)  \right\Vert ^{2}+\left\Vert
U_{j-}^{h};L^{2}\left(  \Pi_{-}\right)  \right\Vert ^{2}  &  =\left\Vert
\chi_{+}u_{j}^{h};L^{2}\left(  \Omega^{h}\right)  \right\Vert ^{2}+\left\Vert
\chi_{-}u_{j}^{h};L^{2}\left(  \Omega^{h}\right)  \right\Vert ^{2}=\\
&  =1+\left\Vert \left(  1-\chi_{+}-\chi_{-}\right)  u_{j}^{h};L^{2}\left(
\Omega^{h}\right)  \right\Vert ^{2}\geq\\
&  \geq1-c\exp\left(  -\left(  3h\right)  ^{-1}\tau\right)  \left\Vert
\mathcal{R}_{\tau}u_{j}^{h};L^{2}\left(  \Omega^{h}\right)  \right\Vert
^{2}\geq1-C\exp\left(  -\left(  3h\right)  ^{-1}\tau\right)  .
\end{align*}
Here $\chi_{\pm}\left(  z\right)  =\chi\left(  1\mp z\right)  $ and
$1-\chi_{+}-\chi_{-}\neq0$ only in $\left(  -\frac{2}{3},\frac{2}{3}\right)
\ni z$ where \ $\mathcal{R}_{\tau}\left(  x\right)  \geq c_{\tau}\exp\left(
\left(  3h\right)  ^{-1}\tau\right)  ,$ $c_{\tau}>0,$ in accord with
(\ref{4.5}). Thus, for a small $h>0,$ we have%
\begin{equation}
\left\Vert U_{j+}^{h};L^{2}\left(  \Pi_{+}\right)  \right\Vert ^{2}+\left\Vert
U_{j-}^{h};L^{2}\left(  \Pi_{-}\right)  \right\Vert ^{2}\geq\frac{1}{2}.
\label{4.12}%
\end{equation}
By (\ref{4.4}), (\ref{4.11}) and (\ref{1.5}), (\ref{1.7}), we obtain%
\[
h^{2}\lambda_{j}^{h}\leq c_{j},
\]%
\begin{align}
\left\Vert U_{j\pm}^{h};H^{1}\left(  \Pi_{\pm}\right)  \right\Vert ^{2}  &
=h^{3}\int_{\Pi_{\pm}}\left(  \left\vert \nabla_{\xi}\left(  \chi u_{j}%
^{h}\right)  \right\vert ^{2}+\left\vert \chi u_{j}^{h}\right\vert
^{2}\right)  d\xi=\label{4.13}\\
&  =\int_{\Omega^{h}}\left(  h^{2}\left\vert \nabla_{x}\left(  \chi_{\pm}%
u_{j}^{h}\right)  \right\vert +\left\vert \chi_{\pm}u_{j}^{h}\right\vert
^{2}\right)  dx\leq\nonumber\\
&  \leq c\left(  h^{2}\left\Vert \nabla_{x}u_{j}^{h};L^{2}\left(  \Omega
^{h}\right)  \right\Vert ^{2}+\left\Vert u_{j}^{h};L^{2}\left(  \Omega
^{h}\right)  \right\Vert ^{2}\right)  =c\left(  h^{2}\lambda_{j}^{h}+1\right)
\leq C_{j}.\nonumber
\end{align}
Thus, the following convergence occurs along an infinitesimal positive
sequence $\left\{  h_{k}\right\}  _{k\in\mathbb{N}}:$%
\begin{equation}
h^{2}\lambda_{j}^{h}\rightarrow\Lambda_{j}^{0}\text{ \ and \ }U_{j\pm}%
^{h}\rightarrow U_{j\pm}^{0}\text{ weakly in }\mathring{H}^{1}\left(  \Pi
_{\pm};\Sigma_{\pm}\right)  . \label{4.14}%
\end{equation}
Unfortunately, we cannot derive from (\ref{4.14}) the strong convergence in
$L^{2}\left(  \Pi_{\pm}\right)  $ in the unbounded domains $\Pi_{\pm}.$
However, we again make use of the weighted estimate (\ref{4.10}) and write%
\begin{align*}
\left\Vert \exp\left(  \tau\zeta^{\pm}\right)  U_{j\pm}^{h};L^{2}\left(
\Pi_{\pm}\right)  \right\Vert  &  \leq\left\Vert \exp\left(  \tau
h^{-1}\left(  1\mp z\right)  \right)  \chi_{\pm}u_{j}^{h};L^{2}\left(
\Omega^{h}\right)  \right\Vert \leq\\
&  \leq c\left\Vert \mathcal{R}_{\tau}u_{j}^{h};L^{2}\left(  \Omega
^{h}\right)  \right\Vert \leq C.
\end{align*}
We now use the compact embedding $H^{1}\left(  \Pi_{\pm}\left(  T\right)
\right)  \subset L^{2}\left(  \Pi_{\pm}\left(  T\right)  \right)  $ in the
finite cylinder $\left\{  \xi^{\pm}\in\Pi_{\pm}:\zeta^{\pm}<T\right\}  $ and
the estimate%
\begin{align*}
\left\Vert U_{j\pm}^{h};L^{2}\left(  \Pi_{\pm}\diagdown\Pi_{\pm}\left(
T\right)  \right)  \right\Vert  &  \leq\exp\left(  -\tau T\right)  \left\Vert
\exp\left(  \tau\zeta^{\pm}\right)  U_{j\pm}^{h};L^{2}\left(  \Pi_{\pm
}\diagdown\Pi_{\pm}\left(  T\right)  \right)  \right\Vert \leq\\
&  \leq C\exp\left(  -\tau T\right)
\end{align*}
with an infinitesimal bound as $T\rightarrow+\infty.$ These yield%
\begin{equation}
U_{j\pm}^{h}\rightarrow U_{j\pm}^{0}\text{ \ strongly in }L^{2}\left(
\Pi_{\pm}\right)  . \label{4.15}%
\end{equation}

We are in position to derive a variational problem for the tripple $\left\{
\Lambda_{j}^{0},U_{j+}^{0},U_{j-}^{0}\right\}  .$ With any $V_{\pm}%
\in\mathring{H}^{1}\left(  \Pi_{\pm};\Sigma_{\pm}\right)  ,$ we take the test
function%
\[
\Omega^{h}\ni x\mapsto v_{\pm}\left(  x\right)  =h^{-3/2}\chi_{\pm}\left(
z\right)  V_{\pm}\left(  h^{-1}y,h^{-1}\left(  1\mp z\right)  \right)
\]
in the integral identity (\ref{1.5}). Multiplying the identity with $h^{2}$
and going over to the stretched coordinates (\ref{3.1}) lead to%
\begin{align}
0  &  =h^{1/2}\left(  \nabla_{x}u_{j}^{h},\nabla_{x}\left(  \chi_{\pm}V_{\pm
}\right)  \right)  _{\Omega^{h}}-h^{1/2}\lambda_{j}^{h}\left(  \chi_{\pm}%
u_{j}^{h},V_{\pm}\right)  _{\Omega^{h}}=\label{4.16}\\
&  =\left(  \nabla_{\xi}U_{j\pm}^{h},\nabla_{\xi}V_{\pm}\right)  _{\Pi_{\pm}%
}-h^{2}\lambda_{j}^{h}\left(  U_{j\pm}^{h},V_{\pm}\right)  _{\Pi_{\pm}%
}+h^{1/2}\left(  \nabla_{x}u_{j}^{h},V_{\pm}\nabla_{x}\chi_{\pm}\right)
_{\Omega^{h}}-h^{1/2}\lambda_{j}^{h}\left(  u_{j}^{h}\nabla_{x}\chi_{\pm
},\nabla_{x}V_{\pm}\right)  _{\Omega^{h}}.\nonumber
\end{align}
The first and second term on the right, by virtue of (\ref{4.14}),
(\ref{4.15}), converge to $\left(  \nabla_{\xi}U_{j\pm}^{0},\nabla_{\xi}%
V_{\pm}\right)  _{\Pi_{\pm}}$ and $\Lambda_{j}^{0}\left(  U_{j\pm}^{h},V_{\pm
}\right)  _{\Pi_{\pm}},$ respectively. Since $\mathcal{R}\left(  x\right)
\geq\exp\left(  \left(  3h\right)  ^{-1}\tau-c_{H}\right)  $ for $x\in
$supp$\left\vert \nabla_{x}\chi_{\pm}\right\vert ,$ the modulo of the third
and fourth terms does not exceed%
\[
c\exp\left(  -\left(  3h\right)  ^{-1}\tau\right)  \left(  h^{2}\left\Vert
\mathcal{R}_{\tau}\nabla_{x}u_{j}^{h};L^{2}\left(  \Omega^{h}\right)
\right\Vert ^{2}+\left\Vert \mathcal{R}_{\tau}u_{j}^{h};L^{2}\left(
\Omega^{h}\right)  \right\Vert ^{2}\right)  \left\Vert V_{\pm};H^{1}\left(
\Pi_{\pm}\right)  \right\Vert .
\]
Hence, the limit passage $k\rightarrow+\infty,$ $h_{k}\rightarrow0^{+}$
converts (\ref{4.16}) into the couple of integral identities (\ref{3.99}).%
\[
\left(  \nabla_{\xi}U_{j\pm}^{0},\nabla_{\xi}V_{\pm}\right)  _{\Pi_{\pm}%
}=\Lambda_{j}^{0}\left(  U_{j\pm}^{0},V_{\pm}\right)  _{\Pi_{\pm}},\ V_{\pm
}\in\mathring{H}^{1}\left(  \Pi_{\pm};\Sigma_{\pm}\right)  ,
\]
while, in view of (\ref{4.12}), at least one of the functions $U_{j\pm}^{0}$
is nontrivial. Since (\ref{4.4}) provides the inequality $\Lambda_{j}^{0}%
\leq\mu_{1}-\tau_{0}^{2},$ we conclude that $\Lambda_{j}^{0}$ coincides with
one of eigenvalues in (\ref{4.1}) in the case $\delta_{1}\geq\tau_{0}^{2}.$

Let us formulate the obtained result.

\begin{proposition}
\label{prop4.3}If the eigenvalue $\lambda_{j}^{h}$ of the problem (\ref{1.5})
meets the condition (\ref{4.4}), the convergence (\ref{4.14}), (\ref{4.15})
occurs while at least one couple $\left\{  \Lambda_{j}^{0},U_{j\pm}%
^{0}\right\}  $ implies an eigenpair of problem (\ref{3.3})-(\ref{3.5}) in the
semi-infinite cylinder $\Pi_{\pm}$ with the variational formulation
(\ref{3.99}).
\end{proposition}

\subsection{Approximation of eigenvalues and eigenfunctions}

Let $\mathcal{H}$ denote the Hilbert space $\mathring{H}^{1}\left(  \Pi_{\pm
};\Sigma_{\pm}\right)  $ equipped with the scalar product%
\begin{equation}
\left\langle u,v\right\rangle =\left(  \nabla_{x}u,\nabla_{x}v\right)
_{\Omega^{h}}+h^{-2}\left(  u,v\right)  _{\Omega^{h}}. \label{5.1}%
\end{equation}
We introduce the positive self-adjoint compact operator $\mathcal{K}$ in
$\mathcal{H}$ by the formula%
\begin{equation}
\left\langle \mathcal{K}u,v\right\rangle =\left(  u,v\right)  _{\Omega^{h}}
\label{5.2}%
\end{equation}
together with the new spectral parameter%
\begin{equation}
\mu=\left(  \lambda+h^{-2}\right)  ^{-1}. \label{5.3}%
\end{equation}
Comparing (\ref{5.1})-(\ref{5.3}) with (\ref{1.5}), we see that the
variational formulation of the problem (\ref{1.3})-(\ref{1.4}) is equivalent
to the abstract equation%
\begin{equation}
\mathcal{K}u=\mu u,\ u\in\mathcal{H}. \label{5.4}%
\end{equation}
The operator $\mathcal{K}$\ has the essential spectrum $\left\{
\mu=0\right\}  $ (cf. \cite[Thm. 9.2.1]{BiSo}) and the positive infinitesimal
sequence $\left\{  \mu_{j}^{h}=\left(  \lambda_{j}^{h}+h^{-2}\right)
^{-1}\right\}  _{j\in N}$ of eigenvalues, while $\widehat{u}_{j}%
^{h}=\left\Vert u_{j}^{h};\mathcal{H}\right\Vert ^{-1}u_{j}^{h}$ are the
corresponding normalized eigenfunctions.

The following assertion is known as the lemma on "almost eigenvalues and
eigenfunctions" (see \cite{ViLu} and, e.g., \cite{BiSo}).

\begin{lemma}
\label{lem5.1}Let $\phi\in\mathcal{H}$ and $\sigma\in\mathbb{R}_{+}$ satisfy
the conditions%
\begin{equation}
\left\Vert \phi;\mathcal{H}\right\Vert =1,\ \ \left\Vert \mathcal{K}%
\phi-\sigma\phi;\mathcal{H}\right\Vert =\theta\in\left(  0,\sigma\right)  .
\label{5.5}%
\end{equation}
Then $\mathcal{K}$ has an eigenvalue in the segment $\left[  \sigma
-\theta,\sigma+\theta\right]  \subset\mathbb{R}_{+}.$ Moreover, for any
$\theta_{1}\in\left(  \theta,\sigma\right)  ,$ one finds coefficients
$a_{k}^{h},...,a_{k+K-1}^{h}$ such that%
\begin{equation}
\left\Vert \phi-\sum_{j=k}^{k+K-1}a_{k}^{h}\widehat{u}_{k}^{h};\mathcal{H}%
\right\Vert \leq2\frac{\theta}{\theta_{1}},\ \ \sum_{j=k}^{k+K-1}\left\vert
a_{k}^{h}\right\vert ^{2}=1, \label{5.6}%
\end{equation}
where $\mu_{k}^{h},...,\mu_{k+K-1}^{h}$ is the complete list of eigenvalues of
$\mathcal{K}$ in the segment $\left[  \sigma-\theta_{1},\sigma+\theta
_{1}\right]  \subset\mathbb{R}_{+}.$
\end{lemma}

The following couples ought to be chosen as approximate solutions of the
spectral equation (\ref{5.2}):%
\begin{equation}
\sigma^{h}=h^{-2}\left(  \Lambda_{j}^{\pm}+1\right)  ^{-1},\ \ \phi
^{h}=\left\Vert v_{j\pm}^{h};\mathcal{H}\right\Vert ^{-1}v_{j\pm}^{h},
\label{5.7}%
\end{equation}
where $\Lambda_{j}^{\pm}$ is an eigenvalue in (\ref{4.1}) with the
corresponding eigenfunction $U_{j}^{\pm}$ of problem (\ref{3.99}) in $\Pi
_{\pm},$ subject to the conditions (\ref{4.312}), and%
\begin{equation}
v_{j\pm}^{h}\left(  x\right)  =h^{-1/2}\chi_{\pm}\left(  z\right)  U_{j}^{\pm
}\left(  \varepsilon^{-1}y,\varepsilon^{-1}\left(  1\mp z\right)  \right)  .
\label{5.8}%
\end{equation}
To estimate the discrepancy $\theta$ of the couple (\ref{5.7}) in the equation
(\ref{5.2}), we need to study the decay properties of the eigenfunction
$U_{j}^{\pm}.$ This can be made on base of the general theory of elliptic
problems in domains with cylindrical outlets to infinity (see the papers
\cite{Ko, Pazy, MaPl2}, \cite{na262} and, e.g., monographs \cite{NaPl,
KoMaRo}). However, for the reader convenience we here present an elementary
proof which is rather similar to our proof of Lemma \ref{lem4.2}.

\begin{lemma}
\label{lem5.2}Let $\left\{  \Lambda^{\pm},U^{\pm}\right\}  $ be an eigenpair
of the problem (\ref{3.99}) in $\Pi_{\pm}$ such that $\left\Vert U^{\pm}%
;L^{2}\left(  \Pi_{\pm}\right)  \right\Vert =1$ and the condition (\ref{4.1})
is met. Then the inclusion $\exp\left(  \tau\zeta^{\pm}\right)  U^{\pm}%
\in\mathring{H}^{1}\left(  \Pi_{\pm};\Sigma_{\pm}\right)  $ and the estimate%
\begin{equation}
\left\Vert \exp\left(  \tau\zeta^{\pm}\right)  U^{\pm};H^{1}\left(  \Pi_{\pm
}\right)  \right\Vert \leq c\left(  \mu_{1}-\Lambda^{\pm}-\tau^{2}\right)
^{-1/2} \label{5.9}%
\end{equation}
are valid with any $\tau\in\left(  0,\left(  \mu_{1}-\Lambda^{\pm}\right)
^{1/2}\right)  .$
\end{lemma}

\textbf{Proof.} We omit the index $\pm.$ We introduce the weight function%
\[
R_{\tau,T}\left(  \xi\right)  =\left\{
\begin{array}
[c]{rr}%
\exp\left(  \tau c_{H}\right)  , & \zeta<c_{H},\\
\exp\left(  \tau\zeta\right)  , & \zeta\in\left[  c_{H},T\right]  ,\\
\exp\left(  \tau T\right)  , & \zeta>T,
\end{array}
\right.
\]
which is continuous and bounded. Moreover, $\nabla_{x}R_{\tau,T}=0$ for
$\zeta\notin\left[  c_{H},T\right]  $ and%
\begin{equation}
\left\vert \nabla_{\xi}R_{\tau,T}\left(  \xi\right)  \right\vert \leq\tau
R_{\tau,T}\left(  \xi\right)  . \label{5.10}%
\end{equation}
Inserting the test function $V=R_{\tau,T}^{2}U\in\mathring{H}^{1}\left(
\Pi;\Sigma\right)  $ into (\ref{3.99}) and repeating the calculation
(\ref{4.7}) with an evident modification, we arrive at the relation%
\begin{align}
\Lambda\left\Vert R_{\tau,T}U;L^{2}\left(  \Pi\right)  \right\Vert ^{2}  &
=\left\Vert \nabla_{\xi}\left(  R_{\tau,T}U\right)  ;L^{2}\left(  \Pi\right)
\right\Vert ^{2}-\left\Vert U\nabla_{\xi}R_{\tau,T};L^{2}\left(  \Pi\right)
\right\Vert ^{2}\geq\label{5.11}\\
&  \geq\mu_{1}\left\Vert R_{\tau,T}U;L^{2}\left(  \Pi_{H}\right)  \right\Vert
^{2}-\tau^{2}\left\Vert R_{\tau,T}U;L^{2}\left(  \Pi_{H}\right)  \right\Vert
^{2},\nonumber
\end{align}
where $\Pi_{H}=\omega\times\left(  c_{H},+\infty\right)  $ and $c_{H}%
=\max\left\{  \left\vert H\left(  \eta\right)  \right\vert :\eta\in
\overline{\omega}\right\}  .$ In (\ref{5.11}) we have applied the Friedrichs
inequality in $\omega$ integrated over $\left(  c_{H},+\infty\right)  \ni
\zeta$ (cf. (\ref{4.8})) and the relation (\ref{5.10}). We finally obtain%
\[
\Lambda\left\Vert U;L^{2}\left(  \Pi\diagdown\Pi_{H}\right)  \right\Vert
^{2}\geq\left(  \mu_{1}-\Lambda-\tau^{2}\right)  \left\Vert R_{\tau,T}%
U;L^{2}\left(  \Pi_{H}\right)  \right\Vert ^{2}%
\]
or, by the normalization assumption on $U,$
\[
\left\Vert R_{\tau,T}U;L^{2}\left(  \Pi_{H}\right)  \right\Vert ^{2}%
\leq\left(  \mu_{1}-\Lambda-\tau^{2}\right)  ^{-1}\Lambda\left\Vert
U;L^{2}\left(  \Pi\diagdown\Pi_{H}\right)  \right\Vert ^{2}\leq\left(  \mu
_{1}-\Lambda-\tau^{2}\right)  ^{-1}\Lambda.
\]
Since the function $\left[  c_{H},+\infty\right)  \ni\zeta\mapsto R_{\tau
,T}\left(  \eta,\zeta\right)  $ is monotone, the limit passage $T\rightarrow
+\infty$ furnishes the inequality%
\[
\left\Vert \exp\left(  \tau\zeta\right)  U;L^{2}\left(  \Pi_{H}\right)
\right\Vert ^{2}\leq\left(  \mu_{1}-\Lambda-\tau^{2}\right)  ^{-1}\Lambda,
\]
which together with the relations%
\begin{align}
\left\Vert R_{\tau,T}\nabla_{\xi}U;L^{2}\left(  \Pi\right)  \right\Vert ^{2}
&  \leq c\left\Vert R_{\tau,T}U;L^{2}\left(  \Pi\right)  \right\Vert
^{2},\label{5.12}\\
\left\Vert \exp\left(  \tau\zeta\right)  U;H^{1}\left(  \Pi\diagdown\Pi
_{H}\right)  \right\Vert ^{2}  &  \leq\left\Vert U;H^{1}\left(  \Pi
\diagdown\Pi_{H}\right)  \right\Vert ^{2}\leq c\left(  \Lambda+1\right)
,\nonumber
\end{align}
inherited from (\ref{5.11}) and (\ref{3.99}), lead to (\ref{5.9}).
$\blacksquare$

To apply Lemma \ref{lem5.1}, we need some calculations. First, using
(\ref{5.9}), we observe that%
\begin{align}
\left\Vert v_{j}^{\pm};\mathcal{H}\right\Vert ^{2}  &  =\frac{1}{h}%
\int_{\Omega^{h}}\chi_{\pm}^{2}\left(  \left\vert \nabla_{x}U_{j}^{\pm
}\right\vert ^{2}+h^{-2}\left\vert U_{j}^{\pm}\right\vert ^{2}\right)
dx+\frac{1}{h}\int_{\Omega^{h}}\nabla_{x}\chi_{\pm}\left(  2U_{j}^{\pm}%
\nabla_{x}U_{j}^{\pm}+\left\vert U_{j}^{\pm}\right\vert ^{2}\nabla_{x}%
\chi_{\pm}\right)  dx\geq\label{5.13}\\
&  \geq\frac{1}{2h}\int_{\Omega^{h}}\chi_{\pm}^{2}\left(  \left\vert
\nabla_{x}U_{j}^{\pm}\right\vert ^{2}+h^{-2}\left\vert U_{j}^{\pm}\right\vert
^{2}\right)  dx-\frac{c}{h}\int_{\Omega^{h}}\left\vert \nabla_{x}\chi_{\pm
}\right\vert \left(  \left\vert \nabla_{x}U_{j}^{\pm}\right\vert
^{2}+\left\vert U_{j}^{\pm}\right\vert ^{2}\right)  dx.\nonumber
\end{align}
By the change of coordinates $x\mapsto\xi^{\pm},$ the first integral on the
right reduces to%
\begin{equation}
\frac{1}{2}\int_{\Pi_{\pm}}\left(  \left\vert \nabla_{\xi}U_{j}^{\pm
}\right\vert ^{2}+\left\vert U_{j}^{\pm}\right\vert ^{2}\right)  d\xi-\frac
{1}{2h}\int_{\Omega^{h}}\left(  1-\chi_{\pm}^{2}\right)  \left(  \left\vert
\nabla_{x}U_{j}^{\pm}\right\vert ^{2}+h^{-2}\left\vert U_{j}^{\pm}\right\vert
^{2}\right)  dx. \label{5.14}%
\end{equation}
The functions $\left\vert \nabla_{x}\chi_{\pm}\right\vert $ and $1-\chi_{\pm
}^{2}$ vanish as $\pm z>\frac{2}{3}$ according to the definition of the
cut-off functions $\chi_{\pm}$ and, therefore, bringing the weight
$\exp\left(  \tau\zeta^{\pm}\right)  \geq\exp\left(  \left(  3h\right)
^{-1}\tau\right)  $ in the last integrals in (\ref{5.13}) and (\ref{5.14})
yields the following upper bound for the integrals:%
\begin{align*}
&  c\exp\left(  \left(  3h\right)  ^{-1}\tau\right)  \frac{1}{h}\int_{\left\{
x\in\Omega^{h}:\pm z>2/3\right\}  }\exp\left(  \tau\zeta^{\pm}\right)  \left(
\left\vert \nabla_{x}U_{j}^{\pm}\right\vert ^{2}+h^{-2}\left\vert U_{j}^{\pm
}\right\vert ^{2}\right)  dx\\
&  \leq c\exp\left(  \left(  3h\right)  ^{-1}\tau\right)  \int_{\Pi_{\pm}}%
\exp\left(  \tau\zeta^{\pm}\right)  \left(  \left\vert \nabla_{\xi}U_{j}^{\pm
}\left(  \xi\right)  \right\vert ^{2}+\left\vert U_{j}^{\pm}\left(
\xi\right)  \right\vert ^{2}\right)  d\xi\leq C\exp\left(  \left(  3h\right)
^{-1}\tau\right)  .
\end{align*}
The exponent $\tau>0$ is taken from Lemma \ref{lem5.2}. By $\left\Vert
U_{j}^{\pm};L^{2}\left(  \Pi_{\pm}\right)  \right\Vert =1$ and (\ref{3.99}),
the first integral in (\ref{5.14}) turns into $\frac{1}{2}\left(
1+\Lambda_{j}^{\pm}\right)  .$ Thus, for a small $h>0,$ we have%
\begin{equation}
\left\Vert v_{j}^{\pm};\mathcal{H}\right\Vert \geq\frac{1}{2}\left(
1+\Lambda_{j}^{\pm}\right)  ^{1/2}. \label{5.15}%
\end{equation}

We further observe that%
\begin{align}
\theta &  =\left\Vert \mathcal{K}\phi-\sigma\phi;\mathcal{H}\right\Vert
=\sup\left\vert \left\langle \mathcal{K}\phi-\sigma\phi;\Psi\right\rangle
\right\vert =\label{5.16}\\
&  =\sigma\left\Vert v_{j}^{\pm};\mathcal{H}\right\Vert ^{-1}\sup\left\vert
h^{-2}\Lambda_{j}^{\pm}\left(  v_{j}^{\pm},\Psi\right)  _{\Omega^{h}}-\left(
\nabla_{x}v_{j}^{\pm},\nabla_{x}\Psi\right)  _{\Omega^{h}}\right\vert
,\nonumber
\end{align}
where supremum is computed over all functions $\Psi\in\mathcal{H}$ such that
$\left\Vert \Psi;\mathcal{H}\right\Vert =1.$ Recalling (\ref{5.8}) and
(\ref{3.1}), we obtain \
\begin{align}
\left(  \nabla_{x}v_{j}^{\pm},\nabla_{x}\Psi\right)  _{\Omega^{h}}%
-h^{-2}\Lambda_{j}^{\pm}\left(  v_{j}^{\pm},\Psi\right)  _{\Omega^{h}}  &
=h^{-1/2}h^{-1}\left\{  \left(  \nabla_{\xi}U_{j}^{\pm},\nabla_{\xi}\left(
\chi_{\pm}\Psi\right)  \right)  _{\Pi_{\pm}}-\Lambda_{j}^{\pm}\left(
U_{j}^{\pm},\chi_{\pm}\Psi\right)  _{\Pi_{\pm}}\right\}  +\label{5.17}\\
&  +h^{-1/2}\left(  U_{j}^{\pm}\nabla_{x}\chi_{\pm},\Psi\right)  _{\Omega^{h}%
}-h^{-1/2}\left(  \nabla_{x}U_{j}^{\pm},\Psi\nabla_{x}\chi_{\pm}\right)
_{\Omega^{h}}.\nonumber
\end{align}
The expression in the curly brackets vanishes by virtue of the integral
identity (\ref{3.99}) with the test function $\xi\mapsto\chi_{\pm}\left(
z\right)  \Psi\left(  x\right)  $ which has a compact support and, therefore,
falls into the function space $\mathring{H}^{1}\left(  \Pi_{\pm};\Sigma_{\pm
}\right)  .$ Modulo of the last two terms in (\ref{5.17}) does not exceed%
\begin{align}
&  ch^{-1/2}\left(  \left\Vert U_{j}^{\pm};L^{2}\left(  \omega\times\left[
-\tfrac{2}{3},\tfrac{2}{3}\right]  \right)  \right\Vert \left\Vert \nabla
_{x}\Psi;L^{2}\left(  \Omega^{h}\right)  \right\Vert +h\left\Vert \nabla
_{x}U_{j}^{\pm};L^{2}\left(  \omega\times\left[  -\tfrac{2}{3},\tfrac{2}%
{3}\right]  \right)  \right\Vert h^{-1}\left\Vert \Psi;L^{2}\left(  \Omega
^{h}\right)  \right\Vert \right) \label{5.177}\\
&  \leq ch\exp\left(  -\left(  3h\right)  ^{-1}\tau\right)  \left\Vert
U_{j}^{\pm};H^{1}\left(  \Pi_{\pm}\right)  \right\Vert \left\Vert
\Psi;\mathcal{H}\right\Vert .\nonumber
\end{align}
Thus, omitting the factor $h$ which is unimportant since $\tau<\left(  \mu
_{1}-\Lambda_{j}^{\pm}\right)  ^{1/2}$ is arbitrary, we conclude the estimate%
\begin{equation}
\theta\leq c\exp\left(  -\left(  3h\right)  ^{-1}\tau\right)  . \label{5.18}%
\end{equation}

Owing to the relationship (\ref{5.3}) for the spectral parameters, Lemma
\ref{lem5.1} extracts an eigenvalue $\lambda_{q}^{l}$ from the sequence
(\ref{1.6}) such that%
\begin{equation}
\left\vert \left(  \lambda_{q}^{h}+h^{-2}\right)  ^{-1}-h^{2}\left(
\Lambda_{j}^{\pm}+1\right)  ^{-1}\right\vert \leq c\exp\left(  -\left(
3h\right)  ^{-1}\tau\right)  . \label{5.19}%
\end{equation}
This inequality transforms into the following one:%
\begin{equation}
\left\vert \lambda_{q}^{h}-h^{-2}\Lambda_{j}^{\pm}\right\vert \leq ch^{-4}%
\exp\left(  -\left(  3h\right)  ^{-1}\tau\right)  \left(  \Lambda_{j}^{\pm
}+1\right)  \left(  h^{2}\lambda_{q}^{h}+1\right)  \leq Ch^{-4}\exp\left(
-\left(  3h\right)  ^{-1}\tau\right)  \label{5.20}%
\end{equation}
where the numbers $C$ and $\tau>0$ depend on $\Lambda_{j}^{\pm}$ only. Here we
have used that $\Lambda_{j}^{\pm}<\mu_{1}$ (cf. (\ref{4.1})) and%
\begin{equation}
\lambda_{q}^{h}\leq h^{-2}\Lambda_{j}^{\pm}+ch^{-4}\exp\left(  -\left(
3h\right)  ^{-1}\tau\right)  \left(  \mu_{1}+1\right)  \left(  h^{2}%
\lambda_{q}^{h}+1\right)  \Rightarrow\lambda_{q}^{h}\leq ch^{-2}\left(
\Lambda_{j}^{\pm}+1\right)  \text{ for a small }h>0. \label{5.21}%
\end{equation}

We now assume that%
\begin{equation}
\Lambda^{\left(  p-1\right)  }<\Lambda^{\left(  p\right)  }=...=\Lambda
^{\left(  p+\varkappa-1\right)  }<\Lambda^{\left(  p+\varkappa\right)  }
\label{5.22}%
\end{equation}
in the family (\ref{4.3}), i.e., $\Lambda^{\left(  p\right)  }$ is an
eigenvalue with multiplicity $\varkappa$. Then the formulas (\ref{5.7}) with
the eigenvalue $\Lambda^{\left(  p\right)  }$ and the corresponding
$\varkappa$ eigenfunctions of the problems (\ref{3.3})-(\ref{3.5}) in
$\Pi_{\pm}$ deliver linear independent approximate solutions $\left\{
h^{-2}\left(  \Lambda^{\left(  p\right)  }+1\right)  ^{-1},\phi^{\left(
m\right)  }\right\}  ,$ $m=p,...,p+\varkappa-1,$ for the spectral equation
(\ref{5.4}). Furthermore, repeating the calculation (\ref{5.13}) and
(\ref{5.17}) with obvious modifications yields the relations%
\begin{equation}
\left\langle \phi^{\left(  m\right)  },\phi^{\left(  k\right)  }\right\rangle
=\delta_{m,k}+O\left(  \exp\left(  -\left(  3h\right)  ^{-1}\tau\right)
\right)  ,\ m,k=p,...,p+\varkappa-1. \label{5.222}%
\end{equation}

We now set $\theta_{1}=\beta\theta$ in the second assertion of Lemma
\ref{lem5.1} where $\beta>1$ is big but will be fixed independent of $h.$ The
condition $\theta_{1}<h^{-2}\left(  \Lambda^{\left(  p\right)  }+1\right)
^{-1}$ in Lemma \ref{lem5.1} can be achieved by diminishing $h.$ As a result,
we obtain coefficients $a_{l}^{hm}$ such that
\begin{equation}
\left\Vert \phi^{\left(  m\right)  }-\sum_{l=k}^{k+K-1}a_{l}^{hm}\widehat
{u}_{l}^{h};\mathcal{H}\right\Vert \leq2\beta^{-1},\ \ \sum_{l=k}%
^{k+K-1}\left\vert a_{l}^{hm}\right\vert ^{2}=1, \label{5.23}%
\end{equation}
where $\lambda_{k}^{h},...,\lambda_{k+K-1}^{h}$ imply all eigenvalues of
problem (\ref{1.5}) such that%
\begin{equation}
\left\vert \left(  \lambda_{l}^{h}+h^{-2}\right)  ^{-1}-h^{-2}\left(
\Lambda^{\left(  p\right)  }+1\right)  ^{-1}\right\vert \leq c\beta\exp\left(
-\left(  3h\right)  ^{-1}\tau\right)  . \label{5.24}%
\end{equation}
We emphasize that, in comparison with Lemma \ref{lem5.1}, we have enlarged in
(\ref{5.24}) the bound $\theta_{1}=\beta\theta$ according to the previous
estimate (\ref{5.18}). However, this does not influence the conclusion: we
chose null coefficients $a_{l}^{hm}$ for new eigenvalues involved.

Recall that $\left\langle \widehat{u}_{l}^{h},\widehat{u}_{j}^{h}\right\rangle
=\delta_{l,j}.$ By (\ref{5.23}) and (\ref{5.22}), we have%
\begin{align}
\left\vert \sum_{l=k}^{k+K-1}a_{l}^{hm}a_{l}^{hr}-\delta_{m,r}\right\vert  &
=\left\vert \left\langle \sum_{l=k}^{k+K-1}a_{l}^{hm}\widehat{u}_{l}^{h}%
,\sum_{j=k}^{k+K-1}a_{j}^{hr}\widehat{u}_{j}^{h}-\delta_{m,r}\right\rangle
\right\vert \leq\left\vert \left\langle \phi^{\left(  m\right)  }%
,\phi^{\left(  r\right)  }\right\rangle -\delta_{m,r}\right\vert
+\label{5.2555}\\
&  +\left\vert \left\langle \phi^{\left(  m\right)  }-\sum_{l=k}^{k+K-1}%
a_{l}^{hm}\widehat{u}_{l}^{h},\sum_{j=k}^{k+K-1}a_{j}^{hr}\widehat{u}_{j}%
^{h}\right\rangle -\delta_{m,r}\right\vert +\left\vert \left\langle
\phi^{\left(  m\right)  },\phi^{\left(  r\right)  }-\sum_{j=k}^{k+K-1}%
a_{j}^{hr}\widehat{u}_{j}^{h}\right\rangle \right\vert \leq\nonumber\\
&  \leq c\exp\left(  -\left(  3h\right)  ^{-1}\tau\right)  +4\beta
^{-1}.\nonumber
\end{align}
Hence, for a small $h>0$ and a large $\beta>1,$ the coefficient columns
$a^{hp},...,a^{hp+\varkappa-1}\in\mathbb{R}^{k}$ are "almost bi-orthogonal"
that can happen in the case $k\geq\varkappa$ only. As a result, we detect at
least $\varkappa$ eigenvalues in (\ref{1.6}) subject to the estimate
(\ref{5.24}). We fix an appropriate $\beta$ and recall that, according to
(\ref{5.21}), the estimate (\ref{5.24}) provides the inequality%
\begin{equation}
\left\vert \left(  \lambda_{l}^{h}-h^{-2}\Lambda^{\left(  p\right)  }\right)
\right\vert \leq ch^{-4}\exp\left(  -\left(  3h\right)  ^{-1}\tau\right)  \leq
c_{p}\exp\left(  -h^{-1}\tau_{p}\right)  \label{5.25}%
\end{equation}
with any $\tau_{p}<\dfrac{\tau}{3}$ and a certain $c_{p}.$ Note that we have
replaced $\dfrac{\tau}{3}$ by $\tau_{p}$ in order to subdue the factor
$h^{-4}$ while Lemma \ref{lem5.2} permits to take any $\tau<\left(  \mu
_{1}-\Lambda^{\left(  p\right)  }\right)  ^{1/2}.$

\begin{proposition}
\label{prop5.7}If the formula (\ref{5.22}) is valid with entries of the
eigenvalue (\ref{4.8}), then at least $\varkappa$ eigenvalues $\lambda_{l}%
^{h}$ of the problem (\ref{1.5}) satisfy the inequality (\ref{5.25}).
\end{proposition}

\subsection{The proof of Theorem \ref{theor4.1}}

It remains to check up that, under the assumption (\ref{5.22}), the
eigenvalues $\lambda_{p}^{h},...,\lambda_{p+\varkappa-1}^{h}$ and no other in
(\ref{1.6}) meet the estimate (\ref{5.25}), i.e. $k=p$ and $K=\varkappa.$ We
know that $K\geq\varkappa.$ Assuming $K>\varkappa,$ we infer the convergence
$h^{2}\lambda_{j}^{h}\rightarrow\Lambda^{\left(  p\right)  }$ for
$j=k,...,k+K-1$ while, by Proposition \ref{prop4.3} and, in particular,
formula (\ref{4.15}), the limits $U_{j\pm}^{0}\in\mathring{H}^{1}\left(
\Pi_{\pm};\Sigma_{\pm}\right)  $ satisfy the problems (\ref{3.99}) in
$\Pi_{\pm}$ and inherit the linear independence from $u_{k}^{h},...,u_{k+K-1}%
^{h}.$ Since $\Lambda^{\left(  p\right)  }$ is an eigenvalue in $\Pi_{\pm}$ of
multiplicity $\varkappa_{\pm}$ and $\varkappa=\varkappa_{+}+\varkappa_{-},$
the above inference and assumption are invalid. In other words, $K=\varkappa.$

If $k>p,$ then, by Proposition \ref{prop5.7}, the total multiplicity of the
spectrum (\ref{1.6}) in the segment%
\begin{equation}
\left[  0,\Lambda^{\left(  p\right)  }+c_{p}\exp\left(  -h^{-1}\tau
_{p}\right)  \right]  \label{5.30}%
\end{equation}
is bigger than $p+\varkappa-1$ for a small $h>0$ and again we readily find out
a contradiction with Proposition \ref{prop4.3}. Finally, the case $k<p$ is
impossible due to Proposition \ref{prop5.7} which, dealing with all
eigenvalues in (\ref{4.3}), detect at least $p+\varkappa-1$ eigenvalues
$\lambda_{l}^{h}$ in the segment (\ref{5.30}). Thus, $k=p$ and Theorem
\ref{theor4.1} is proved.

\subsection{Asymptotic expansions for eigenfunctions and the localization
effect\label{sect3.5}}

Let $\Lambda^{\left(  p\right)  }$ be an eigenvalue of multiplicity
$\varkappa$ (see (\ref{5.22})) in the list (\ref{4.3}) and let $p+\varkappa
-1<N.$ Then by Theorem \ref{theor4.1} there exists $h_{p}>0$ such that for
$h\in\left(  0,h_{p}\right]  $ the segment%
\begin{equation}
\left[  \frac{1}{2}\left(  \Lambda^{\left(  p\right)  }+\Lambda^{\left(
p-1\right)  }\right)  ,\frac{1}{2}\left(  \Lambda^{\left(  p\right)  }%
+\Lambda^{\left(  p+\varkappa\right)  }\right)  \right]  \label{5.31}%
\end{equation}
contains the eigenvalues $\lambda_{p}^{h},...,\lambda_{p+\varkappa-1}^{h}$ of
the problem (\ref{1.5}) and is free of other entries in the eigenvalue
sequence (\ref{1.6}). Owing to the relationship (\ref{5.3}), we see that the
eigenvalues $\mu_{p}^{h}=\left(  \lambda_{p}^{h}+h^{-2}\right)  ^{-1}%
,...,\mu_{p}^{h}=\left(  \lambda_{p+\varkappa-1}^{h}+h^{-2}\right)  ^{-1}$ and
no other eigenvalue of the operator $\mathcal{K}$ fall into the segment%
\begin{equation}
\left[  2h^{2}\left(  \Lambda^{\left(  p\right)  }+\Lambda^{\left(
p+\varkappa\right)  }+2\right)  ^{-1},2h^{2}\left(  \Lambda^{\left(  p\right)
}+\Lambda^{\left(  p-1\right)  }+2\right)  ^{-1}\right]  . \label{5.32}%
\end{equation}
We now apply the second assertion in Lemma \ref{lem5.1} with $\theta_{1}%
=h^{2}\theta_{1}^{0}$ where $\theta_{1}^{0}$ is independent of $h$ and,
moreover, the segment $\left[  h^{2}\left(  \Lambda^{\left(  p\right)
}+1\right)  -h^{2}\theta_{1}^{0},h^{2}\left(  \Lambda^{\left(  p\right)
}+1\right)  +h^{2}\theta_{1}^{0}\right]  $ lies inside (\ref{5.32}). As a
result we find coefficients $a_{l}^{hr}$ such that, according to (\ref{5.7})
and (\ref{5.15}),%
\begin{equation}
\left\Vert \phi^{\left(  l\right)  }-\sum_{r=p}^{p+\varkappa-1}a_{l}%
^{hr}\widehat{u}_{r}^{h};\mathcal{H}\right\Vert \leq2\frac{\theta}{\theta_{1}%
}\leq c_{p}h^{-2}\exp\left(  -\left(  3h\right)  ^{-1}\tau\right)
,\ l=p,...,p+\varkappa-1. \label{5.33}%
\end{equation}
We compare formulas (\ref{5.1}), (\ref{5.4}) and (\ref{1.5}), (\ref{1.7}) to
derive that%
\[
\left\Vert \widehat{u}_{l}^{h};\mathcal{H}\right\Vert ^{2}=\left(  \lambda
_{l}^{h}+h^{-2}\right)  \left\Vert u_{j}^{h};L^{2}\left(  \Omega^{h}\right)
\right\Vert ^{2}=\lambda_{l}^{h}+h^{-2},\ \ l=p,...,p+\varkappa-1.
\]
Furthermore, by repeating the calculations (\ref{5.13}) and (\ref{5.17})
similarly to (\ref{5.177}), we obtain that%
\[
\left\vert \left\Vert v_{l}^{h};\mathcal{H}\right\Vert ^{2}-\left(
1+\Lambda^{\left(  l\right)  }\right)  \right\vert \leq c\exp\left(  -\left(
3h\right)  ^{-1}\tau\right)
\]
where $v_{l}^{h}$ denotes the approximate solution (\ref{5.8}) corresponding
to the eigenvalue $\Lambda^{\left(  l\right)  }=\Lambda_{j}^{\pm}$ (see
(\ref{5.7}) and cf. (\ref{4.1}), (\ref{4.3})). Since $\Lambda_{j}^{\pm}$ and
$\lambda_{p}^{h},...,\lambda_{p+\varkappa-1}^{h}$ satisfy (\ref{5.19}), the
inequalities (\ref{5.33}) convert into%
\begin{equation}
\left\Vert v_{l}^{h}-h^{2}\sum_{r=p}^{p+\varkappa-1}\widehat{a}_{l}^{hr}%
u_{r}^{h};\mathcal{H}\right\Vert \leq c_{p}\exp\left(  -\left(  3h\right)
^{-1}\tau\right)  \label{5.34}%
\end{equation}
while, by means of the formula (\ref{5.2555}) with $\beta=\theta_{1}/\theta,$
the estimates%
\begin{equation}
\left\vert \sum_{q=p}^{p+\varkappa-1}\widehat{a}_{l}^{hq}\widehat{a}_{r}%
^{hq}-\delta_{l,r}\right\vert \leq c_{p}h^{-2}\exp\left(  -\left(  3h\right)
^{-1}\tau\right)  ,\ \ l,r=p,...,p+\varkappa-1, \label{5.35}%
\end{equation}
are valid.

Now we employ the following simple algebraic fact (see, e.g.,
\cite[Lemma 7.1.7]{Nabook}, \cite[Lemma 3.3]{na317}): under the
condition (\ref{5.35}) one can find out an orthogonal matrix
$b^{h}=\left(  b_{l}^{hq}\right)  $ such
that the inequalities (\ref{5.34}) ensure the estimates%
\begin{equation}
\left\Vert u_{l}^{h}-h^{-2}\sum_{r=p}^{p+\varkappa-1}b_{l}^{hr}v_{r}%
^{h};\mathcal{H}\right\Vert \leq c_{p}h^{-2}\exp\left(  -\left(  3h\right)
^{-1}\tau\right)  ,\ l=p,...,p+\varkappa-1. \label{9.1}%
\end{equation}

\begin{theorem}
\label{theor3.987}If $\Lambda^{\left(  p\right)  }$ is an eigenvalue of
multiplicity $\varkappa$ (see (\ref{4.3}) and (\ref{5.22})), then there exists
the coefficient columns $b_{l}^{h}=\left(  b_{l}^{hp},...,b_{l}^{hp+\varkappa
-1}\right)  ,\ l=p,...,p+\varkappa-1,$ composing an orthogonal matrix of size
$\varkappa\times\varkappa$ and furnishing the inequalities (\ref{9.1}) where
$u_{p}^{h},...,u_{p+\varkappa-1}^{h}$ are eigenfunctions of problem
(\ref{1.5}) under the normalization and orthogonality conditions (\ref{1.7})
and $v_{p}^{h},...,v_{p+\varkappa-1}^{h}$ are the functions (\ref{5.8})
constructed from eigenfunctions of problem (\ref{3.99}) in the semi-infinite
cylinders $\Pi_{\pm}$ (see (1.25)) under the normalization and orthogonality
conditions (\ref{4.312}).
\end{theorem}

Since the eigenfunctions $U_{j}^{\pm}$ decay exponentially in the
semi-cylinders $\Pi_{\pm}$ as $\zeta^{\pm}\rightarrow+\infty$ (see Lemma
\ref{lem5.2}), the functions $v_{p}^{h}=v_{j\pm}^{h}$ in (\ref{9.1}) and
(\ref{5.8}) are of order $h^{-1/2}$ in the vicinity of the end $\Gamma
_{h}^{\pm}$ of the thin cylinder $\Omega^{h}$ but become exponentially small
at a distance from this end. Thus, the estimate (\ref{9.1}) obtained in
Theorem \ref{theor3.987} exhibits the localization effect discussed in Section
\ref{sect1.2}. On the other hand, the structure of the weight function
(\ref{4.5}), which is exponentially large in the middle of the thin cylinder
(\ref{1.1}), and the estimate (\ref{4.10}) of weighted norms of an
eigenfunction $u^{h}$ of problem (\ref{1.5}) ensure the same effect. We
emphasize that in the case when $\Lambda^{\left(  p\right)  }=\Lambda_{j}%
^{\pm}$ is a simple eigenvalue in the list (\ref{4.3}), Theorem
\ref{theor3.987} provides the localization of the corresponding eigenfunction
$u_{p}^{h}$ in a neighborhood of only one end $\Gamma_{h}^{\pm}$ while Lemma
\ref{lem4.2} cannot distinguish between the ends $\Gamma_{h}^{+}$ and
$\Gamma_{h}^{-}$ (see (\ref{1.2})). In Section \ref{sect4.1} we demonstrate an
example of the principal eigenfunction $u_{1}^{h}$ which does not become
exponentially small near the both ends.

\section{Similar spectral problems\label{sect4}}

\subsection{Splitting of a multiple eigenvalue in the case $n=2$%
\label{sect4.1}}

If $\omega=\left(  0,1\right)  \subset\mathbb{R}^{1}$ and
\begin{equation}
\Omega^{h}=\left\{  x=\left(  y,z\right)  \in\mathbb{R}^{2}:\eta=h^{-1}%
y\in\left(  0,1\right)  ,\ z<\left(  -1-hH_{-}\left(  \eta\right)
,1+hH_{+}\left(  \eta\right)  \right)  \right\}  , \label{7.1}%
\end{equation}
then the first eigenpair of the Dirichlet problem (\ref{1.9}) in the interval
$\left(  0,1\right)  $ becomes%
\[
\mu_{1}=\pi^{2},\ \varphi_{1}\left(  \eta\right)  =2^{-1/2}\sin\left(  \pi
\eta\right)  .
\]
Hence, the condition (\ref{3.13}) converts into%
\begin{equation}
0>\int_{0}^{1}H_{\pm}\left(  \eta\right)  \left(  \cos^{2}\left(  \pi
\eta\right)  -\sin^{2}\left(  \pi\eta\right)  \right)  d\eta=\int_{0}%
^{1}H_{\pm}\left(  \eta\right)  \cos\left(  2\pi\eta\right)  d\eta.
\label{7.2}%
\end{equation}
By Theorem \ref{theor3.1}, the two-dimensional problem (\ref{3.3})-(\ref{3.5})
gets a trapped mode and the localization effect occurs in the case when the
coefficient (\ref{7.2}) in the Fourier series at least of one function
$H_{\pm}\in C\left[  0,1\right]  $ is negative.

For the symmetric thin domain (\ref{7.1}) , i.e., for $H_{+}=H_{-},$ the
simple first eigenvalue $\Lambda_{1}=\Lambda_{1}^{\pm}\in\sigma_{d}\left(
A_{\pm}\right)  $ of the operator $A_{\pm}$ of the problem (\ref{3.3}%
)-(\ref{3.5}) in $\Pi=\Pi_{\pm}$ must be treated as multiple in Theorem
\ref{theor4.1} which gives the same asymptotic forms for $\lambda_{1}^{h}$ and
$\lambda_{2}^{h}$ in (\ref{1.6}) due to the equality $\Lambda^{\left(
1\right)  }=\Lambda^{\left(  2\right)  }$ in (\ref{4.3}). At the same time,
the first eigenvalue $\lambda_{1}^{h}$ of the problem (\ref{1.3})-(\ref{1.4})
is simple by the maximum principle. To display a difference in the asymptotics
of $\lambda_{1}^{h}$ and $\lambda_{2}^{h},$ we need to construct the second
term in the expansion of the eigenvalues.

The Fourier method ensures that, for $\zeta>c_{H},$ the eigenfunction $U_{1}$
of the problem (\ref{3.3})-(\ref{3.5}) admits the decomposition
\begin{equation}
U_{1}\left(  \xi\right)  =C_{1}\exp\left(  -\left(  \pi^{2}-\Lambda
_{1}\right)  ^{1/2}\zeta\right)  \sin\left(  \pi\eta\right)  +O\left(
\exp\left(  -\left(  4\pi^{2}-\Lambda_{1}\right)  ^{1/2}\zeta\right)  \right)
. \label{7.3}%
\end{equation}
The factor $C_{1}$ in (\ref{7.3}) does not vanish because in the case
$C_{1}=0$ the next term
\[
C_{2}\exp\left(  -\left(  4\pi^{2}-\Lambda_{1}\right)  ^{1/2}\zeta\right)
\sin\left(  2\pi\eta\right)
\]
in the decomposition of $U_{1}$ becomes main but changes sign.

Following \cite[\S 5.6]{MaNaPl}, we accept the asymptotic ans\"{a}tze for an
eigenpair $\left\{  \lambda^{h},u^{h}\right\}  $ of the problem (\ref{1.3}%
)-(\ref{1.4})%
\begin{align}
\lambda^{h}  &  =h^{-2}\left(  \Lambda_{1}+\epsilon\Lambda_{1}^{\prime
}+...\right)  ,\label{7.4}\\
u^{h}  &  =\sum_{\pm}\left(  b_{\pm}U_{1}\left(  h^{-1}y,h^{-1}\left(  1\mp
z\right)  \right)  \right)  +\epsilon U_{\pm}^{\prime}\left(  h^{-1}%
y,h^{-1}\left(  1\mp z\right)  +...\right) \nonumber
\end{align}
where the dots stand for neglectible terms and%
\begin{equation}
\epsilon=\exp\left(  -2h^{-1}\left(  \pi^{2}-\Lambda_{1}\right)  ^{1/2}%
\zeta\right)  . \label{7.5}%
\end{equation}
We insert (\ref{7.4}) into the equation (\ref{1.3}) and the Neumann boundary
conditions (\ref{1.4}) while remarking that the Dirichlet conditions
(\ref{1.34}) are satisfied. In view of (\ref{3.1}) we write%
\begin{align*}
\exp\left(  -\left(  \pi^{2}-\Lambda_{1}\right)  ^{1/2}\zeta^{\pm}\right)   &
=\exp\left(  -h^{-1}\left(  \pi^{2}-\Lambda_{1}\right)  ^{1/2}\left(  1\mp
z\right)  \right)  =\exp\left(  -h^{-1}\left(  \pi^{2}-\Lambda_{1}\right)
^{1/2}\left(  2-\left(  1\mp z\right)  \right)  \right)  =\\
&  =\exp\left(  -2h^{-1}\left(  \pi^{2}-\Lambda_{1}\right)  ^{1/2}\right)
\exp\left(  \left(  \pi^{2}-\Lambda_{1}\right)  ^{1/2}\zeta^{\mp}\right)
\end{align*}
that explains the form (\ref{7.3}) of the new small parameter. Furthermore, we
derive the following problems in $\Pi$ to determine the second terms in the
asymptotic ans\"{a}tze:%
\begin{align}
-\Delta_{\xi}U_{\pm}^{\prime}\left(  \xi\right)  -\Lambda_{1}U_{\pm}^{\prime
}\left(  \xi\right)   &  =b_{\pm}\Lambda_{\pm}^{\prime}U_{1}\left(
\xi\right)  ,\ \xi\in\Pi,\label{7.6}\\
U_{\pm}^{\prime}\left(  \xi\right)   &  =0,\ \xi\in\Sigma,\nonumber\\
\partial_{\nu}U_{\pm}^{\prime}\left(  \xi\right)   &  =-b_{\mp}\partial_{\nu
}\left(  \exp\left(  \left(  \pi^{2}-\Lambda_{1}\right)  ^{1/2}\zeta^{\pm
}\right)  \sin\left(  \pi\eta\right)  \right)  ,\ \xi\in\Gamma.\nonumber
\end{align}
Since the eigenvalues $\Lambda_{\pm}^{\prime}$ are simple, the only
compatibility condition in the problem (\ref{7.6}) reads:%
\begin{equation}
\Lambda_{\pm}^{\prime}b_{\pm}=\Lambda_{\pm}^{\prime}b_{\pm}\int_{\Pi}%
U_{1}\left(  \xi\right)  ^{2}d\xi=b_{\mp}\int_{\Pi}U_{1}\left(  \xi\right)
\partial_{\nu}\left(  \exp\left(  \left(  \pi^{2}-\Lambda_{1}\right)
^{1/2}\zeta^{\pm}\right)  \sin\left(  \pi\eta\right)  \right)  d\eta=:b_{\mp
}F. \label{7.7}%
\end{equation}
The authors do not know a way to confirm the assumption
\begin{equation}
F\neq0. \label{7.8}%
\end{equation}
However, under this assumption, the system of the two ($\pm$) linear algebraic
equations (\ref{7.7}) has the following couple of solutions%
\[
\Lambda_{\pm}^{\prime}=-F,\ b_{\pm}=\pm2^{-1/2}\text{ and }\Lambda_{\pm
}^{\prime}=F,\ b_{\pm}=2^{-1/2}%
\]
that complete the asymptotic ans\"{a}tze (\ref{7.4}). A
straightforward and simple modification of the arguments in
Section \ref{sect3} justifies the
constructed asymptotics and, in particular, gives the relations%
\begin{equation}
\lambda_{1}^{h}=h^{-2}\left(  \Lambda_{1}-\epsilon F+O\left(  \epsilon
^{3/2}\right)  \right)  \text{ \ and \ }\lambda_{2}^{h}=h^{-2}\left(
\Lambda_{1}+\epsilon F+O\left(  \epsilon^{3/2}\right)  \right)  . \label{7.9}%
\end{equation}

\begin{figure}
[ptb]
\begin{center}
\includegraphics[
height=1.356in,
width=3.1592in
]%
{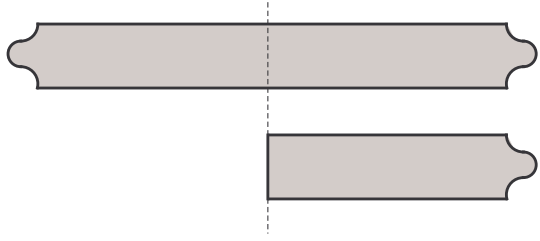}%
\caption{Division of a symmetric domain}%
\label{f8}%
\end{center}
\end{figure}

The formulas (\ref{7.9}), (\ref{7.8}) show the asymptotic fission of the first
two entries in the eigenvalue sequence (\ref{1.6}) of the problem
(\ref{1.3})-(\ref{1.4}) in the symmetric thin domain $\Omega^{h}.$ We note
that $\lambda_{1}^{h}$ and $\lambda_{2}^{h},$ respectively, are eigenvalues of
the spectral problem restricted on the half domain (see Fig. \ref{f8})%
\[
\Omega^{h+}=\left\{  x=\left(  y,z\right)  :y\in\left(  0,h\right)
,\ 0<z<1+H\left(  h^{-1}y\right)  \right\}
\]
with the Neumann and Dirichlet conditions on the straight lateral side while,
respectively, the even and odd in $z$ extensions of the eigenfunctions
complete eigenpairs of the problem (\ref{1.5}) in the domain (\ref{7.1}). In
this way one sees that the localization effect occurs at the both ends
simultaneously. However, this observation does not help to verify the
inequality (\ref{7.8}).%

\subsection{The Dirichlet problem\label{sect4.2}}

Let us consider the equation (\ref{1.3}) with the Dirichlet boundary
conditions on the whole boundary
\begin{equation}
u^{h}\left(  x\right)  =0,\ x\in\partial\Omega^{h}. \label{7.10}%
\end{equation}
By the classical approach \cite{Rellich}, the Dirichlet problem%
\begin{equation}
-\Delta_{\xi}U=\Lambda U\text{ \ in \ }\Pi_{\pm},\ \ U=0\text{\ on \ }%
\partial\Pi_{\pm}, \label{7.11}%
\end{equation}
in the semi-infinite cylinder
\begin{equation}
\Pi_{\pm}=\left\{  \xi^{\pm}:\eta\in\omega,\ \zeta^{\pm}>-H_{\pm}\left(
\eta\right)  \right\}  \label{7.12}%
\end{equation}
with the profile function $H_{\pm}\in C^{1}\left(  \overline{\omega}\right)
$, has no solution in $\mathring{H}^{1}\left(  \Pi_{\pm};\partial\Pi_{\pm
}\right)  .$ In other words, the unbounded operator $A_{\pm}^{D},$ generated
by the quadratic form (\ref{3.6}) in $\mathring{H}^{1}\left(  \Pi_{\pm
};\partial\Pi_{\pm}\right)  $ (see \cite[\S 10.2]{BiSo}), possesses empty
point spectrum. As a result, we see that the localization effect discovered in
the previous sections for the mixed boundary value problem (\ref{1.3}%
)-(\ref{1.4}) under the hypotheses of Theorems \ref{theor3.1} and
\ref{theor3.2} does not occur in the same domain (\ref{1.1}). However,
employing a result in \cite{Jones}, we observe the effect in the case of a
dumbbell domain in Fig. \ref{f3}.

\begin{figure}
[ptb]
\begin{center}
\includegraphics[
height=1.0525in,
width=3.0969in
]%
{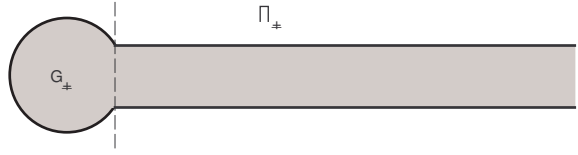}%
\caption{Infinite cane-head domain}%
\label{f9}%
\end{center}
\end{figure}

Let $G_{\pm}\subset\mathbb{R}_{-}^{n}=\left\{  \xi^{\pm}:\zeta^{\pm
}<0\right\}  $ be a domain such that the boundary $\partial G_{\pm}$ contains
the set $\left\{  \xi:\zeta^{\pm}=0,\ \eta\in\omega\right\}  $ and the closure
$\overline{G_{\pm}}$ is compact. We set
\begin{align}
\Pi_{\pm}  &  =G_{\pm}\cup\left(  \omega\times\left[  0,+\infty\right)
\right)  ,\label{7.13}\\
\Omega^{h}  &  =\left(  \omega\times\left[  -1,1\right]  \right)  \cup%
{\displaystyle\bigcup_{\pm}}
G_{\pm}^{h},\ G_{\pm}^{h}=\left\{  x:\xi^{\pm}\in G_{\pm}\right\}
\label{7.14}%
\end{align}
(see Fig. \ref{f9} and \ref{f3}, respectively).%

We assume that the first eigenvalue $\Lambda_{1\pm}^{D}$ of the Dirichlet
problem in the bounded domain $G_{\pm}$ lies below the cut-off $\mu_{1}$ of
the problem (\ref{7.11}), (\ref{7.10}) in the unbounded domain (\ref{7.12}).
One may fulfill this requirement by inflating a domain of a fixed shape.
Following \cite{Jones}, we extend the first eigenfunction $\Phi_{1}^{D\pm}$ by
zero from $G_{\pm}$ onto $\Pi_{\pm}$ and apply the minimum principle
(\cite[Thm. 10.2.1]{BiSo}) to derive that%
\begin{align}
\min\left\{  \Lambda:\Lambda\in\sigma\left(  A_{\pm}^{D}\right)  \right\}   &
=\inf_{U\in\mathring{H}^{1}\left(  \Pi_{\pm};\partial\Pi_{\pm}\right)  }%
\frac{a\left(  U,U\right)  }{\left\Vert U;L^{2}\left(  \Pi_{\pm}\right)
\right\Vert ^{2}}\leq\frac{\left\Vert \nabla_{\xi}\phi_{1}^{D\pm};L^{2}\left(
\Pi_{\pm}\right)  \right\Vert ^{2}}{\left\Vert \phi_{1}^{D\pm};L^{2}\left(
\Pi_{\pm}\right)  \right\Vert ^{2}}=\label{7.15}\\
&  =\frac{\left\Vert \nabla_{\xi}\phi_{1}^{D\pm};L^{2}\left(  G_{\pm}\right)
\right\Vert ^{2}}{\left\Vert \phi_{1}^{D\pm};L^{2}\left(  G_{\pm}\right)
\right\Vert ^{2}}=\Lambda_{1\pm}^{D}<\mu_{1}.\nonumber
\end{align}

As in Section \ref{sect2}, the relation (\ref{7.15}) assures that the discrete
spectrum of operator $A_{\pm}^{D}$ is not empty and the problem (\ref{7.11})
in the cane-head domain (\ref{7.13}) has an eigenvalue $\Lambda_{1}\in\left(
0,\mu_{1}\right)  $. In this way, by the inflation of $G_{\pm}$, one can place
any given number of eigenvalues in the interval $\left(  0,\mu_{1}\right)  $.

We repeat word by word the argumentation in Section \ref{sect3} and arrive at
the following assertion.

\begin{theorem}
\label{theor7.9} Let the operators $A_{\pm}^{D}$ of the problems (\ref{7.11})
in $\Pi_{\pm}$ have the eigenvalues (\ref{4.1}). Then Theorem \ref{theor4.1}
keeps the validity for the eigenvalue sequence (\ref{1.6}) of the Dirichlet
problem (\ref{1.3}), (\ref{7.10}) in the dumbbell domain $\Omega^{h}$ in
(\ref{7.14}).
\end{theorem}

\subsection{The Neumann problem}

Let us change (\ref{7.10}) for the Neumann boundary condition%
\begin{equation}
\partial_{\nu}u\left(  x\right)  =0,\ \ x\in\partial\Omega^{h}, \label{7.16}%
\end{equation}
where $\partial_{\nu}$ stands for the differentiation along the outward normal
$\nu$ which is defined almost everywhere on the Lipschitz (by our assumption)
boundary of the domain (\ref{1.1}) (see Fig. \ref{f1}).

The integral identity%
\begin{equation}
\left(  \nabla_{x}u^{h},\nabla_{x}v\right)  _{\Omega^{h}}=\lambda^{h}\left(
u^{h},v\right)  _{\Omega^{h}},\ \ v\in H^{1}\left(  \Omega^{h}\right)  ,
\label{7.17}%
\end{equation}
posed an the whole Sobolev space $H^{1}\left(  \Omega^{h}\right)  $, serves
for the Neumann problem (\ref{1.3}), (\ref{7.10}). It is well-known that the
eigenvalues%
\begin{equation}
0=\lambda_{0}^{h}<\lambda_{1}^{h}\leq\lambda_{2}^{h}\leq...\leq\lambda_{p}%
^{h}\leq...\rightarrow+\infty\label{7.18}%
\end{equation}
of the problem (\ref{7.17}) satisfy the asymptotic formula%
\begin{equation}
\left\vert \lambda_{p}^{h}-\frac{\pi^{2}}{4}p^{2}\right\vert \leq
\mathbf{c}h^{\frac{1}{2}}p^{3}\text{ for }h\in\left(  0,hp^{-4}\right)  ,
\label{7.19}%
\end{equation}
where the positive constants $\mathbf{c}$ and $h$ are independent of the
eigenvalue number $p\in%
\mathbb{N}
$ (see \cite[\S 7.1]{Nabook} for the estimation of the bounds in (\ref{7.19})
and also \cite{na355, na370} for similar results in other singular perturbed
spectral problems). Note that $M_{p}=\frac{1}{4}\pi^{2}p^{2},$ $p=0,1,...,$
are eigenvalues of the one-dimensional limit problem%
\begin{equation}
\partial_{z}^{2}w\left(  z\right)  =M_{p}w\left(  z\right)  ,z\in\left(
-1,1\right)  ,\ \partial_{z}w\left(  \pm1\right)  =0, \label{7.20}%
\end{equation}
obtained by the traditional procedure of the dimension reduction (see
\cite{SPSH, Ciar, Nabook} and others).

For $p=O\left(  h^{-\frac{1}{2}}\right)  ,$ formula (\ref{7.19}) no longer
displays an asymptotics of the eigenvalue $\lambda_{p}^{h}$ because the term
$\frac{1}{4}\pi^{2}p^{2}$ becomes of the same order $h^{-1}$ as the bound
$\mathbf{c}h^{\frac{1}{2}}p^{3}$. As mentioned, e.g., in \cite{LoPe, na355},
the above observation can help to reveal other asymptotically stable series of eigenvalues.

To indicate such eigenvalues, we apply the method of artificial boundary
conditions proposed in \cite{Vas1} to find out the point spectrum in the
continuous spectrum of the Neumann Laplacian in a symmetric strip with a
finite obstacle in the inside.

\begin{figure}
[ptb]
\begin{center}
\includegraphics[
height=0.493in,
width=5.444in
]%
{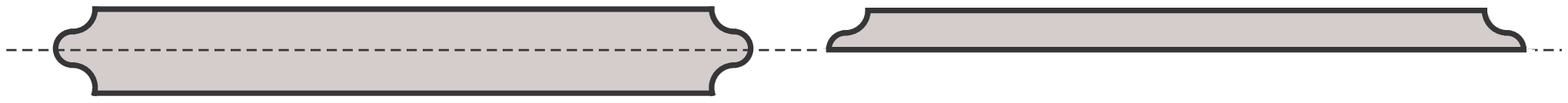}%
\caption{Division of a symmetric domain}%
\label{f10}%
\end{center}
\end{figure}

Let us assume that the thin domain (\ref{1.1}) is symmetric with respect to
the plane $\left\{  x:y_{1}=0\right\}  ,$ i.e., formally we have%
\begin{equation}
\Omega^{h}=\left\{  x:\left(  -y_{1},y_{2},...,y_{n-1},z\right)  \in\Omega
^{h}\right\}  . \label{7.21}%
\end{equation}
Then we restrict the equation (\ref{1.3}) on $\Omega_{\wedge}^{h}=\left\{
x\in\Omega^{h}:y_{1}>0\right\}  $ and the boundary condition (\ref{7.16}) on
$\left(  \partial\Omega^{h}\right)  _{\wedge}=\left\{  x\in\partial\Omega
^{h}:y_{1}>0\right\}  $ (see Fig. \ref{f10}) while imposing the artificial
boundary conditions%
\begin{equation}
u_{\wedge}^{h}=0\text{ on }\Upsilon^{h}=\left\{  x\in\Omega^{h}:y_{1}%
=0\right\}  . \label{7.22}%
\end{equation}
%

The variational formulation of the new mixed boundary value problem refers to
the integral identity%
\begin{equation}
(\nabla_{x}u_{\wedge}^{h},\nabla_{x}v_{\wedge})_{\Omega_{\wedge}^{h}}%
=\lambda_{\wedge}^{h}\left(  u_{\wedge}^{h},v_{\wedge}\right)  _{\Omega
_{\wedge}^{h}},v\in\mathring{H}^{1}\left(  \Omega_{\wedge}^{h};\Upsilon
^{h}\right)  . \label{7.23}%
\end{equation}
Let us make several observation. First, any eigenpair $\left\{  \lambda
_{\wedge}^{h},u_{\wedge}^{h}\right\}  $ of the problem (\ref{7.23}) becomes an
eigenpair of the Neumann problem (\ref{7.17}) after the extension of
$u_{\wedge}^{h}$ from $\Omega_{\wedge}^{h}$ onto $\Omega^{h}$ as an odd
function in $y_{1}.$ Second, in the limit boundary value problem%
\begin{align}
-\Delta_{\xi}U^{\wedge}  &  =\Lambda^{\wedge}U^{\wedge}\text{ \ in \ }%
\Pi^{\wedge}=\left\{  \xi\in\Pi:\eta_{1}>0\right\}  ,\label{7.24}\\
U^{\wedge}  &  =0\text{\ on \ }\Upsilon^{\wedge}=\left\{  \xi\in\Pi:\eta
_{1}=0\right\}  ,\ \ \partial_{\nu}U^{\wedge}=0\text{ \ on \ }\partial
\Pi^{\wedge}\diagdown\overline{\Upsilon^{\wedge}},\nonumber
\end{align}
posed in the half of the cylinder (\ref{3.98}), the continuous spectrum
$\left[  \mu_{1}^{\wedge},+\infty\right)  $ begins with the first eigenvalue
$\mu_{1}^{\wedge}>0$ of the problem%
\begin{align}
-\Delta_{\xi}\varphi^{\wedge}  &  =\mu^{\wedge}\varphi^{\wedge}\text{ \ in
\ }\omega^{\wedge}=\left\{  \eta\in\omega:\eta_{1}>0\right\}  ,\label{7.25}\\
\varphi^{\wedge}  &  =0\text{\ on \ }\upsilon^{\wedge}=\left\{  \eta\in
\omega:\eta_{1}=0\right\}  ,\ \ \partial_{\nu}\varphi^{\wedge}=0\text{ \ on
\ }\partial\omega^{\wedge}\diagdown\overline{\upsilon^{\wedge}},\nonumber
\end{align}
(compare with (\ref{3.3})-(\ref{3.5}) and (\ref{1.10})). Third, owing to the
positive cut-off $\mu_{1}^{\wedge},$ the arguments used in Section
\ref{sect2.3} to prove Theorem \ref{theor3.1} maintain the following assertion.

\begin{theorem}
\label{theor7.23}Let the cross-section $\omega$ be symmetric, i.e.,
$\omega=\left\{  \eta:\left(  -\eta_{1},\eta_{2},...,\eta_{n-1}\right)
\in\omega\right\}  $, and let the even in $\eta_{1}$ function $H\in C\left(
\overline{\omega}\right)  $ meet the condition%
\begin{equation}
\int_{\omega^{\wedge}}H\left(  \eta\right)  \left(  \left\vert \nabla_{\eta
}\varphi_{1}^{\wedge}\left(  \eta\right)  \right\vert ^{2}-\mu_{1}^{\wedge
}\varphi_{1}^{\wedge}\left(  \eta\right)  ^{2}\right)  d\eta<0, \label{7.26}%
\end{equation}
where $\left\langle \mu_{1}^{\wedge},\varphi_{1}^{\wedge}\right\rangle $ is
the first eigenpair of the problem (\ref{7.25}). Then the mixed boundary value
problem (\ref{7.24}) has an eigenvalue below the cut-off $\mu_{1}^{\wedge}$
for the continuous spectrum.
\end{theorem}

\begin{remark}
\label{rem7.24}We have used in (\ref{3.97}) the formula%
\begin{equation}
2\left(  \nabla_{\eta}\varphi_{1},\varphi_{1}\nabla_{\eta}H\right)  _{\omega
}=\left(  \varphi_{1}^{2},\Delta_{\eta}H\right)  _{\omega} \label{7.27}%
\end{equation}
in order to derive the simplified condition (\ref{3.15}) in Theorem
\ref{theor3.2}. In general the change $\varphi_{1},\omega\mapsto\varphi
_{1}^{\wedge},\omega^{\wedge}$ makes (\ref{7.27}) wrong, since $\varphi
_{1}^{\wedge}$ does not vanish on $\upsilon^{\wedge}.$ However, in the case
$\nabla_{\eta}H\left(  \eta\right)  =0$ for $\eta\in\upsilon$ the assertion in
Theorem \ref{theor3.2} remains valid under the symmetry assumption
(\ref{7.21}). $\blacksquare$
\end{remark}

The Dirichlet conditions on the part $\upsilon^{\wedge}$ of the boundary are
sufficient to realize the justification scheme, developed in Section
\ref{sect3}. Indeed, the key inequality (\ref{4.8}) ought to be replaced by
the inequality%
\[
\left\Vert \nabla_{\eta}w;L^{2}\left(  \omega^{\wedge}\right)  \right\Vert
^{2}\geq\mu_{1}^{\wedge}\left\Vert w;L^{2}\left(  \omega^{\wedge}\right)
\right\Vert ^{2},\ w\in\mathring{H}^{1}\left(  \omega^{\wedge};\upsilon
^{\wedge}\right)  .
\]
Hence, the above observations provide the following assertion.

\begin{theorem}
\label{theor7.25}Let $\Lambda_{\pm}\in\left(  0,\mu_{1}^{\wedge}\right)  $ be
an eigenvalue of the problem (\ref{7.24}) in $\Pi_{\pm}^{\wedge}=\left\{
\xi^{\pm}:\eta\in\omega^{\wedge},\ \zeta^{\pm}>-H_{\pm}\left(  \eta\right)
\right\}  $ (cf. Theorem \ref{theor7.23} and Remark \ref{rem7.24}). Then the
Neumann problem (\ref{7.23}) has an eigenvalue $\lambda_{N\left(  h\right)
}^{h}$ such that%
\begin{equation}
\left\vert \lambda_{N\left(  h\right)  }^{h}-h^{-2}\Lambda_{\pm}^{\wedge
}\right\vert \leq c_{\wedge}\exp\left(  -h^{-1}\tau_{\wedge}\right)
\label{7.28}%
\end{equation}
where $c_{\wedge}$ and $\tau_{\wedge}$ are certain positive numbers.
\end{theorem}

We emphasize that the eigenvalue number $N\left(  h\right)  $ in Theorem
\ref{theor7.25} depends on the small parameter $h$ and $N\left(  h\right)
\rightarrow+\infty,$ as $h\rightarrow0^{+}.$ Moreover, we cannot assert that
the only eigenvalue $\lambda_{N\left(  h\right)  }^{h}$ satisfies the relation
(\ref{7.28}). These issues follow from the fact that the first series of
eigenvalues with the stable asymptotics (\ref{7.19}) is of order $h^{0}$ but
the eigenvalue $\lambda_{N\left(  h\right)  }^{h}$ is of order $h^{-2}.$

\end{document}